\documentclass[a4paper,
               10pt,
               pdftex,
               normalheadings,
               headsepline,
               footsepline,
               onecolumn,
               headinclude,
               footinclude,
               DIV14,
               abstracton]
{scrartcl}

\usepackage
{
    graphicx,
    amssymb,
    amsmath,
    amsthm,
    xcolor,
    dsfont,
    algpseudocode,
    authblk,
}

\usepackage[ngerman, english]{babel}

\usepackage[
left=15mm,
right=15mm,
top=15mm,
bottom=30mm
]{geometry}

\usepackage[bf]{caption}
\usepackage[colorlinks,
            pdffitwindow=false,
            plainpages=false,
            pdfpagelabels=true,
            pdfpagemode=UseOutlines,
            pdfpagelayout=SinglePage,
            bookmarks=false,
            colorlinks=true,
            hyperfootnotes=false,
            linkcolor=blue,
            citecolor=green!50!black]
{hyperref}

\usepackage{enumerate}
\usepackage{algorithm}

\usepackage{tabularx}

\usepackage{hyphenat}
\DeclareMathAlphabet{\mathpzc}{OT1}{pzc}{m}{it}

\newcommand{\R}{\mathbb{R}}
\newcommand{\N}{\mathbb{N}}

\newcommand{\Fcal}{\mathcal{F}}
\newcommand{\Kcal}{\mathcal{K}}

\newcommand{\Ycal}{\mathcal{Y}}

\newtheorem{theorem}{Theorem}[section]
\newtheorem{corollary}[theorem]{Corollary}

\newtheorem{remark}[theorem]{Remark}
\newtheorem{example}[theorem]{Example}

\makeatletter
\renewcommand*\env@matrix[1][*\c@MaxMatrixCols c]{%
  \hskip -\arraycolsep
  \let\@ifnextchar\new@ifnextchar
  \array{#1}}
\makeatother

\newcolumntype{Y}{>{\centering\arraybackslash}X}

\let\OLDthebibliography\thebibliography
\renewcommand\thebibliography[1]{
	\OLDthebibliography{#1}
	\setlength{\parskip}{0pt}
	\setlength{\itemsep}{0pt plus 0.3ex}
}

\date{}

\begin{document}

\title{Controlling nonlinear PDEs using low-dimensional bilinear approximations obtained from data}
\author[1]{Sebastian Peitz}
\affil[1]{\normalsize Department of Mathematics, Paderborn University, Germany}

\twocolumn
\maketitle

\paragraph{Abstract ---}
\textit{ \hspace{-0.5cm}
In a recent article, we presented a framework to control nonlinear partial differential equations (PDEs) by means of Koopman operator based reduced models and concepts from switched systems \cite{PK19}. The main idea was to transform a control system into a set of autonomous systems for which the optimal switching sequence has to be computed. These individual systems can be approximated very efficiently by reduced order models obtained from data, and one can guarantee equality of the full and the reduced objective function under certain assumptions. In this article, we extend these results to continuous control inputs using convex combinations of multiple Koopman operators corresponding to constant controls, which results in a bilinear control system. Although equality of the objectives can be carried over when the PDE depends linearly on the control, we show that this approach is also valid in other scenarios using several flow control examples of varying complexity.
}


\section{Introduction}
\label{sec:Introduction}

Real-time control of partial differential equations (PDEs) is a highly challenging task, in particular for nonlinear systems, see, e.g., \cite{BN15} for a recent survey on turbulent flow control. To this end, advanced control techniques such as \emph{Model Predictive Control (MPC)} \cite{GP17} or machine learning based control \cite{DBN17} have gained more and more attention in recent years. In MPC, an open-loop optimal control is computed repeatedly on a finite-time horizon using a model of the system dynamics. This results in a feedback loop, but the open-loop problem has to be solved in a very short time, which is in general infeasible for complex systems such as nonlinear PDEs, at least when using classical discretization approaches such as finite elements.

A remedy to this problem are surrogate models which can be solved significantly faster, see \cite{LMQR14,BGW15} for overviews. Besides classical approaches such as \emph{Proper Orthogonal Decomposition (POD)} \cite{Sir87,KV99,Row05,HV05,BDPV18}, an approach that has attracted a lot of attention in recent years is to construct a \emph{reduced order model (ROM)} by means of the \emph{Koopman operator} \cite{Koo31}. This is an infinite-dimensional but linear operator describing the dynamics of observables, and the approach can even be applied to sensor measurements or in situations where the underlying system dynamics is unknown. Significant advances have been made over the past years both theoretically (see, e.g., \cite{Mez05,BMM12}) as well as numerically. In the latter case, \emph{Dynamic Mode Decomposition (DMD)} \cite{Sch10,RMB+09,TRL+14,KGPS18} or \emph{Extended Dynamic Mode Decomposition (EDMD)} \cite{WKR15,KKS16} are the most popular algorithms. 

More recently, various attempts have been made to use Koopman operator based ROMs for control \cite{PBK15,PBK18,BBPK16,KM18b,KKB17}. In many of these approaches, the Koopman operator is approximated for an augmented state (consisting of the actual state and the control) in order to deal with the non-autonomous control system. For this reason, large amounts of data are necessary to cover a sufficient range of the dynamics.
An alternative approach has been presented in \cite{PK19}, where the control system is replaced by a set of autonomous systems with constant inputs. This way, the optimal control problem is transformed into a switching time problem, and equality of the full and the reduced objective function can be shown by utilizing a recent convergence result for EDMD \cite{KM18a}. The major advantage of this approach is that the required amount of data is very small.
However, a drawback is that the resulting control problem is of combinatorial nature and that the input is restricted to a finite set.

In this article, we extend these results by a straightforward interpolation procedure, which results in a bilinear reduced control system. Instead of switching between the autonomous dynamics, we also allow intermediate control values, which is closely related to relaxation approaches known in mixed-integer optimal control \cite{Sag09}. In the reduced model, these controls are realized by linear interpolation between the different Koopman operators. For control affine systems, we prove convergence in measure for the reduced objective function. Embedded in an MPC framework, we observe that the approach also yields remarkable results for nonlinear control dependencies. As an example, we consider the \emph{fluidic pinball} \cite{DPMN18}, a fluid flow problem governed by the 2D incompressible Navier--Stokes equations which shows chaotic behavior.

The remainder of the article is structured as follows. In Section~\ref{sec:Koopman}, we introduce basic concepts for the Koopman operator and its numerical approximation. The interpolation based control is then introduced in Section~\ref{sec:Continuous} and the application to reduced order modeling of PDEs is discussed in Section~\ref{sec:PDE}. Finally, the combination with MPC is addressed in Section \ref{sec:MPC} before we conclude with a short summary and possible future directions in Section~\ref{sec:Conclusion}.

\section{Koopman Operator based optimal control}
\label{sec:Koopman}

The overall goal we pursue is to efficiently solve optimal control problems constrained by PDEs:
\begin{equation}\label{eq:OCP}
	\begin{aligned}
	\min_{u \in \mathcal{U}} \int_{t_0}^{t_e} &L(y(\cdot,t)) + \alpha \|u(t)\|_2 + \beta \|\dot{u}(t)\|_2 \ dt \\
	\mbox{s.t.} \quad \dot{y}(\cdot,t) &= G(y(\cdot,t),u(t)), \\
	y(\cdot,0) &= y^0,
	\end{aligned}
\end{equation}
where $y$ is the system state and $y(\cdot, t)$ is an element of an appropriate function space $\Ycal$ (e.g., the Sobolev space $H^s(\Omega, \R^{n_y})$ with $\Omega$ being the domain, $n_y$ the spatial dimension and $s \geq 1$ the required differentiability). Furthermore, $u \in L^2([t_0,t_e], U)$ is the control with box constraints $U = [u^l, u^u]$, and $G \colon \Ycal \times U \rightarrow \Ycal$ describes the system dynamics. The parameters $\alpha, \beta \in \R^{\geq0}$ are weights for the penalty terms for the control input as well as its variation. For ease of notation, the objective function $L \colon \Ycal\rightarrow\R$ only depends on $y$ explicitly.

For real systems, it is often insufficient to determine open-loop control inputs. A remedy to this issue is MPC \cite{GP17}, where optimal control problems are solved repeatedly on finite horizons. Using a model of the system dynamics, an open-loop problem is solved in real-time over a so-called \emph{prediction horizon} of length $p$. Following \cite{GP17}, we consider discrete dynamics obtained, for instance, by introducing the time-$T$ map of the PDE:
\begin{equation}\label{eq:PDE_DT}
	\begin{aligned}
		y_{i+1} &= \Phi(y_{i},u_i), \\
		y_0 &= y^0,
	\end{aligned}
\end{equation}
The motivation behind this is that the control is constant over each sample time interval such that it is sufficient to consider the flow map $\Phi$ of the continuous dynamics. This results in the following optimal control problem:
\begin{align}
	\min_{u \in U^p} \sum_{i=s}^{s+p-1} &L(y_i) + \alpha \|u_i\|_2 + \beta \| u_i - u_{i-1} \|_2 \label{eq:MPC} \tag{MPC} \\
	\text{s.t.}\quad y_{i+1} &= \Phi(y_i, u_{i-s+1}) \quad \text{for}~i = s,\ldots,s + p - 1, \notag \\
	y_s &= y^s. \notag
\end{align}
The first part of the solution of \eqref{eq:MPC} is then applied to the real system while the optimization is repeated with the prediction horizon moving forward by one sample time. (The indexing $i-s+1$ is required to account for the finite-horizon control and the infinite-horizon state.)

The main challenge in MPC is that Problem~\eqref{eq:MPC} has to be solved within one sample time $h$ which can already be very challenging for ODE constraints. When the system dynamics is described by a PDE, additional measures have to be taken in order to achieve real-time applicability. One such approach is via the Koopman operator, which will be described next.

\subsection{Koopman operator and EDMD}
Let $ \Phi : \mathcal{M} \to \mathcal{M} $ be a discrete deterministic dynamical system defined on the state space $\mathcal{M}$ 
and let $ f \colon \mathcal{M} \rightarrow \R^q $ be a real-valued observable of the system. Then the Koopman operator $ \Kcal \colon \Fcal \to \Fcal $ (with $\Fcal = L^{2}(\mathcal{M})$) is defined by
\begin{equation*}
(\Kcal f)(y) = f(\Phi(y)),
\end{equation*} 
see~\cite{BMM12,Mez13,WKR15} for more details. 
The Koopman operator is linear but infinite-dimensional. It describes the evolution of the observable $ f $ whereas its adjoint, the Perron--Frobenius operator, describes the evolution of densities. 

The most popular numerical approach to approximate the Koopman operator is Dynamic Mode Decomposition (DMD), which was initially proposed by Peter Schmid \cite{Sch10}. It is a modal decomposition method for large data sets such as fluid flow simulations. While being similar to Proper Orthogonal Decomposition (POD) \cite{Sir87}, the main difference is that instead of obtaining a basis with minimal $L^2$ projection error, each of the DMD modes possesses a frequency with which it oscillates, determined by the corresponding complex eigenvalue \cite{RMB+09}. Consequently, DMD can be interpreted as a generalized Fourier transform.

\begin{figure}[t]
	\centering
	\includegraphics[width=0.45\textwidth]{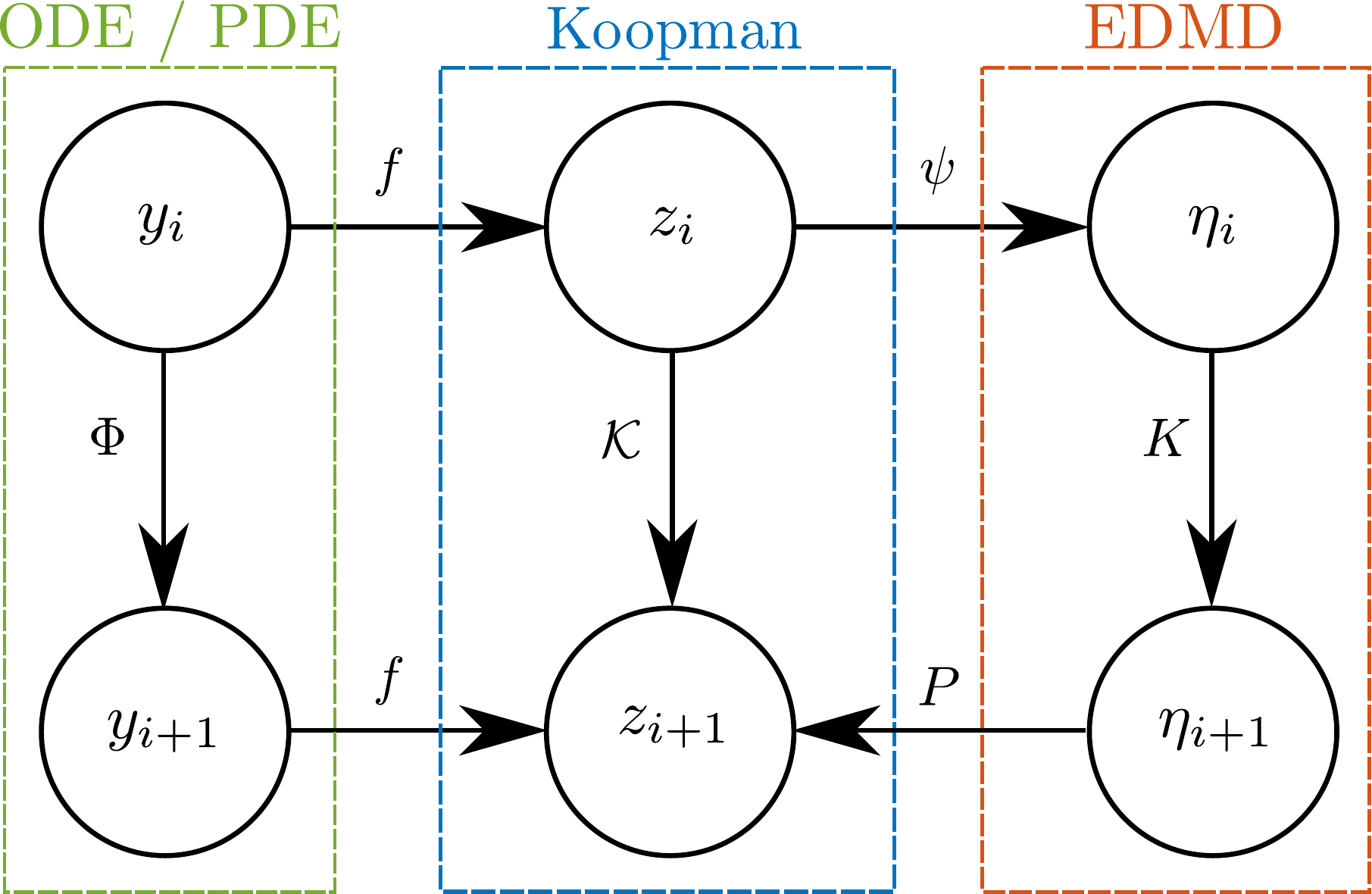}
	\caption{Relation between the system dynamics $\Phi$, the corresponding Koopman operator $\Kcal$ and its finite-dimensional representation $K$ computed via EDMD.}
	\label{fig:Koopman_EDMD}
\end{figure}

\emph{Extended Dynamic Mode Decomposition (EDMD)} \cite{WKR15,KKS16} is a generalization of DMD and can be used to compute a finite-dimensional approximation of the Koopman operator, its eigenvalues, eigenfunctions, and modes, cf.~\cite{KNK+18} for a survey.
EDMD constructs an approximation of the Koopman operator from data (i.e., measurements) given by $ z = f(y) \in \R^q$. For finite-dimensional systems, it is possible to observe the entire (discretized) state (\emph{full state observable}) but the approach is valid for arbitrary observables.
These observations are \emph{lifted} to a higher-dimensional space using arbitrary basis functions (e.g., monomials, Hermite polynomials or radial basis functions). 
For a given set of basis functions $ \{ \psi_{1},\,\psi_{2},\,\dots,\,\psi_{k} \} $ (a so-called \emph{dictionary}), we define $ \psi \colon \R^q \to \R^k $ by
\begin{equation*}
	\psi(z) =
	\begin{bmatrix}
	\psi_{1}(z) & \psi_{2}(z) & \dots & \psi_{k}(z)
	\end{bmatrix}^{\top}.
\end{equation*}
Note that for $ \psi(z) = z $, we recover the standard DMD. The available measurement or simulation data (either obtained from one long or multiple short time series) is written in matrix form as
\begin{equation*}
	Z =
	\begin{bmatrix}
	z_1 & z_2 & \cdots & z_m
	\end{bmatrix}
	\quad\text{and}\quad
	\widetilde{Z} =
	\begin{bmatrix}
	\widetilde{z}_1 & \widetilde{z}_2 & \cdots & \widetilde{z}_m
	\end{bmatrix},
\end{equation*}
where $ \widetilde{z}_i = f(\Phi(y_i)) $. Using the dictionary, the data is embedded into the feature space by
\begin{equation*}
	\begin{aligned}
		\Psi_{Z} &=
		\begin{bmatrix} \psi(z_{1}) & \psi(z_{2}) & \dots & \psi(z_{m}) \end{bmatrix}
		\quad \text{and} \\
		\Psi_{\widetilde{Z}} &=
		\begin{bmatrix} \psi(\widetilde{z}_{1}) & \psi(\widetilde{z}_{2}) & \dots & \psi(\widetilde{z}_{m}) \end{bmatrix},
	\end{aligned}
\end{equation*}
from which we compute the matrix $ K \in \R^{k \times k} $ defined by
\begin{equation*}
K^{\top} = \Psi_{\widetilde{Z}} \Psi_Z^+ = \big( \Psi_{\widetilde{Z}} \Psi_Z^{\top} \big) \big(\Psi_Z \Psi_Z^{\top}\big)^+.
\end{equation*}
This matrix $ K $ can be seen as a finite-dimensional approximation of the Koopman operator. 
The decomposition of the Koopman operator into eigenvalues, eigenfunctions, and modes is commonly used to analyze the system dynamics as well as predict the future state. Alternatively, we can compute updates for the observable $z$ directly using $K$. Introducing $\eta = \psi(z)$, we have
\begin{equation*}
	\eta_{i+1} = K^{\top} \eta_{i}, \quad i = 0,1,\ldots
\end{equation*}
From here, we can obtain $z_{i+1}$ using the projection matrix $P$, cf.~Figure~\ref{fig:Koopman_EDMD}, where the relation between the dynamical system $\Phi$, the related Koopman operator $\Kcal$, and the EDMD approximation $K$ is visualized.

\subsection{Convergence of EDMD towards the Koopman operator}
\label{subsec:Koopman_Convergence}
First results showing convergence of the EDMD algorithm towards the Koopman operator have recently been proved in~\cite{WKR15,AM17,KM18a}. In short, the results state that the matrix $K$ obtained by EDMD converges to the Koopman operator in an $L^2$ sense as both the basis size as well as the number of measurements tend to infinity. To show this, the assumptions have to be satisfied that the basis functions $\psi$ are orthonormal, that the data $z$ is drawn from the space of all measurements $ \mathcal{Z} \subset \R^q $ according to a probability measure $\mu$ such that
\begin{align*}
\mu \{z \in \mathcal{Z} \mid c^{\top}\psi(z) = 0 \} = 0, \quad \mbox{for all}~c \in \R^k, c \neq 0,
\end{align*}
and that the Koopman operator is bounded. Note that in particular the last assumption is hard to verify for many systems.

\subsection{Switched system control of PDEs via K-ROMs}
\label{subsec:Koopman_MPC}
As mentioned at the beginning of this Section, we want to solve Problem~\eqref{eq:MPC} in real-time.
According to the approach in \cite{PK19}, we realize this by
taking two steps:
\begin{enumerate}[i)]
	\item replace the dynamical control system in \eqref{eq:MPC} by a finite number of autonomous systems
	\[\Phi_{\overline{u}}(y_i) = \Phi(y_i,\overline{u})\]
	with constant input $\overline{u} \in \hat{U} = \{u^0,\ldots,u^{n_c-1}\}$;
	\item construct linear systems for low-dimensional observations of the infinite-dimensional systems $\Phi_{u^j}$ using the Koopman operator.
\end{enumerate}
Due to step i), Problem \eqref{eq:MPC} is transformed into a switching problem, where the objective is to select the optimal right hand side in each time step:
\begin{align} 
	\min_{u \in \hat{U}^p} \sum_{i=s}^{s+p-1} &L(y_i) + \alpha \|u_i\|_2 + \beta \| u_i - u_{i-1} \|_2 \label{eq:MPC_STO} \tag{MPCs} \\
	\mbox{s.t.}\quad y_{i+1} &= \Phi_{u_{i-s+1}}(y_i) \quad \text{for}~i = s,\ldots,s + p - 1, \notag \\
	y_s &= y^s. \notag
\end{align}
In other words, each entry of $u$ describes which system $\Phi_{u_i}$ to apply in the $i^{\mathsf{th}}$ step. Problem~\eqref{eq:MPC_STO} is a combinatorial problem that can be solved using dynamic programming \cite{BD15,XA00}, for instance.

Solving either \eqref{eq:MPC} or \eqref{eq:MPC_STO} numerically (e.g., using a finite volume discretization) can quickly become very expensive such that real-time applicability is infeasible. Furthermore, there are many systems where the dynamics is not explicitly known. In both situations, we can use observations (i.e., data) to approximate the Koopman operator and derive a linear system describing the dynamics of these observations. These could consist of (part of) the system state as well as arbitrary functions of the state such as the lift coefficient of an object within a flow field.

Following step ii), we use a \emph{Koopman operator based reduced order model (K-ROM)} to overcome the issue of large computational cost. To this end, we compute $n_c$ Koopman operators for the $n_c$ different autonomous systems:
\begin{equation*}
	(\Kcal_{u^j} f)(y) = f(\Phi_{u^j}(y)),\quad j = 0,\ldots,n_c-1.
\end{equation*}
Using EDMD, we can compute an approximation of the individual Koopman operators from observations of the respective systems and thereby derive linear systems for the observations $\eta = \psi(z) = \psi(f(y))$:
\begin{equation}\label{eq:Discrete_Koopman_Dynamics}
	\eta_{i+1} = {K_{u^j}}^{\top} \eta_{i},\quad j = 0,\ldots,n_c-1.
\end{equation}
These linear dynamics now replace the original time-$T$ map which results in the K-ROM based formulation of \eqref{eq:MPC_STO}:
\begin{align} 
	\min_{u \in \hat{U}^p} \sum_{i=s}^{s+p-1} &\hat{L}(\eta_i) + \alpha \|u_i\|_2 + \beta \| u_i - u_{i-1} \|_2 \label{eq:MPC_STO_Koopman} \tag{K-MPCs}\\
	\mbox{s.t.}\quad \eta_{i+1} &= K_{u_{i-s+1}}^{\top} \eta_{i} \quad \text{for}~i = s,\ldots,s + p - 1, \notag \\
	\eta_s &= \psi(f(y^s)). \notag 
\end{align}
The key advantage over other approaches where one operator is computed for an augmented state $\hat{z} = (z,u)^\top$ (see, e.g., \cite{PBK15,PBK18,KM18b}) is that the individual models \eqref{eq:Discrete_Koopman_Dynamics} can be approximated with very low data requirements. Less than 100 data points for each system are often sufficient.

Using the $L^2$ convergence result for the Koopman operator \cite{KM18a}, we can show convergence in measure for the reduced objective function:
\begin{theorem}[\cite{PK19}]\label{thm:Convergence_MPC}
	Consider Problem \eqref{eq:MPC_STO} and the Koopman operator based approximation \eqref{eq:MPC_STO_Koopman} and assume that we have convergence of EDMD towards the Koopman operator according to  \cite{KM18a}.
	Furthermore, assume $L(y(t)) = \hat{L}(\eta_{i})$ for all $t\in[t_0,t_e]$ and the corresponding $i = (t-t_0)/h$.
	Then, as the basis size $k$ and number of sampled data points $m$ tend to infinity,
	the objective function value of Problem \eqref{eq:MPC_STO_Koopman} converges in measure (with respect to the initial condition $z_0 = f(y^0)$) to that of Problem \eqref{eq:MPC_STO} for every $u \in \hat{U}^p$, $p < \infty$. That is, for all $\epsilon > 0$ we obtain
	\begin{equation*}
	\lim_{k\rightarrow\infty} \lim_{m\rightarrow\infty} \mu \Big( \Big\{ z_0 \in \mathcal{Z} \Big | \sum_{i=0}^{p} | \hat{L}(\eta_i) - L(y_i) | \geq \epsilon \Big\} \Big) = 0.
	\end{equation*}
\end{theorem}
Note that the assumption $L(y(t)) = \hat{L}(\eta_{i})$ is not restrictive in practical settings since we can only consider quantities in the objective function which we can observe. 
By this approach, we can significantly accelerate the computation which is on the one hand due to the linearity of the model and on the other hand due to the restriction to a small number of observables instead of the full state $y$. 

\section{Continuous control inputs}
\label{sec:Continuous}

The switched systems approach presented in Section~\ref{subsec:Koopman_MPC} yields speed-ups of several orders of magnitude. However, two drawbacks are that the resulting optimization problem is of combinatorial nature -- and is thereby harder to solve -- and that the control input is limited to a small number of values, i.e., to the finite set $\hat{U} \subset U$. In order to overcome both of these drawbacks, we use relaxation and define the matrices
\begin{align*}
	A = K^{\top}_{u^0}, \quad B &= \Big[ B_1 \quad \cdots \quad B_{n_c-1} \Big] \\
	&= \Big[ K^{\top}_{u^1} - K^{\top}_{u^0} \quad \cdots \quad K^{\top}_{u^{n_c-1}} - K^{\top}_{u^0} \Big].
\end{align*}
Using $A$ and $B$, we introduce the bilinear control system
\begin{equation} \label{eq:KROM_continuous} \tag{K-ROM}
	\begin{aligned}
	\eta_{i+1} &= A \eta_{i} + \sum_{j=1}^{n_c-1} B_j \eta_{i} \frac{u_{j,i} - u^0}{u^j-u^0}, \\
	\eta_0 &= \psi(f(y^0)),
	\end{aligned}
\end{equation}
where the term bilinear refers to the fact that \eqref{eq:KROM_continuous} contains terms $\eta u_j$, $j=1,\ldots,n_c-1$, but is otherwise linear both in $\eta$ and in $u$ \cite{Ell09}. By relaxing $u_i$ to the convex hull of the individual control inputs, i.e.,
\begin{equation*}
	\resizebox{1\hsize}{!}{$
	\begin{aligned}
		u_{i} \in \overline{U} = \Big\{ \overline{u} \in \R^{n_c-1} ~\Big |~ \overline{u}_{j} \geq 0 \, \forall \, j\in\{1,\ldots, n_c-1\}, \sum_{j=1}^{n_c-1} \overline{u}_{j} \leq 1 \Big\},  
	\end{aligned}
	$}
\end{equation*}
this system simply interpolates linearly between the autonomous dynamics corresponding to the constant inputs $u^0,\ldots,u^{n_c-1}$. Note that \eqref{eq:KROM_continuous} yields the exact dynamics almost everywhere for the individual systems due to the convergence result for EDMD.
For control affine systems, we also have convergence for intermediate values of $u$: 
\begin{theorem}\label{thm:Koopman_continuous}
	Consider a dynamical control system of the form \eqref{eq:PDE_DT} and let $\Kcal_{u^0}, \ldots, \Kcal_{u^{n_c-1}}$ be the Koopman operators associated with the constant inputs $u^0,\ldots,u^{n_c-1}$.
	
	Assume that the observation map $f$ is linear and that the system dynamics $\Phi$ is linear in $u$. Then a convex combination of operators is equal to the Koopman operator for the convex combination of the controls $u^0,\ldots,u^{n_c-1}$, i.e.,
	\begin{equation*}
	\sum_{j=0}^{n_c-1}\gamma_j \Kcal_{u^j} = \Kcal_{\overline{u}},
	\end{equation*}
	where
	\begin{equation*}
	\overline{u} = \sum_{j=0}^{n_c-1}\gamma_j u^j , \quad \sum_{j=1}^{n_c-1}\gamma_j \leq 1, ~ \gamma_j \geq 0 ~ \forall ~ j\in\{0,\ldots, n_c-1\}.
	\end{equation*}
\end{theorem}
\begin{proof}
	The claim follows directly from the linearity:
	\begin{align*}
	&\left(\left(\sum_{j=0}^{n_c-1}\gamma_j \Kcal_{u^j}\right) f\right) y = \sum_{j=0}^{n_c-1}\left(\left(\gamma_j \Kcal_{u^j}\right) f\right) y \\ = &\sum_{j=0}^{n_c-1} \gamma_j f(\Phi(y, u^j)) = f(\sum_{j=0}^{n_c-1} \gamma_j \Phi(y, u^j)) \\ = &f(\Phi(y, \sum_{j=0}^{n_c-1} \gamma_j u^j)) = f(\Phi(y,  \overline{u})) = \big(\Kcal_{\overline{u}} f\big)y.
	\end{align*}
\end{proof}

An important consequence of Theorem \ref{thm:Koopman_continuous} is that we can use \eqref{eq:KROM_continuous} to predict the dynamics of the observations $z$.
\begin{corollary}\label{cor:EqualDynamics}
	Denote by $\varphi: \Ycal \times \R \times {U}^i \rightarrow \Ycal$ and $\hat{\varphi}: \R^k \times \R \times {U}^i \rightarrow \R^k$ the (discrete-time) flows of the full system \eqref{eq:PDE_DT} and of the reduced dynamics \eqref{eq:KROM_continuous}, respectively.
	Then, as $m$ and $k$ tend towards infinity,	the state of the reduced model \eqref{eq:KROM_continuous} converges in measure (with respect to the initial condition $z_0 = f(y^0)$) towards the observations of \eqref{eq:PDE_DT}, i.e.,
	\begin{equation*}
	\resizebox{1\hsize}{!}{$
	\begin{aligned}
	&\lim_{k,m\rightarrow\infty} \mu \left( \left\{ z_0 \in \mathcal{Z} \mid \| \hat{\varphi}(\eta_0, t_i, u) - \psi(f(\varphi(y^0, t_i, u))) \| \geq \epsilon \right\} \right) = 0,
	\end{aligned}
	$}
	\end{equation*}
	for $\epsilon > 0$ and $i\in\{1,\ldots,p\}$, $p< \infty$. 
\end{corollary}
\begin{proof}
	The assumptions yield $L^2$ convergence of EDMD towards the Koopman operator according to \cite{KM18a}, which implies convergence in measure. That is, for all $\epsilon > 0$ we obtain for all $u \in \{u_0, \ldots, u_{n_c-1} \}^p$:
	\begin{equation*}
	\resizebox{1\hsize}{!}{$
	\begin{aligned}
	&\lim_{k\rightarrow\infty}\lim_{m\rightarrow\infty} \mu \left( \left\{ z_i \in \mathcal{Z} \mid \|K_{u_i}^{\top} \psi(z_i) - \mathcal{K}_{u_i} \psi(z_i) \| \geq \epsilon \right\} \right) \\ 
	=&\lim_{k\rightarrow\infty}\lim_{m\rightarrow\infty} \mu \left( \left\{ z_i \in \mathcal{Z} \mid \|K_{u_i}^{\top} \psi(z_i) - \psi(f(\Phi_{u_i}(y(\cdot,t_{i})))) \| \geq \epsilon \right\} \right) \\ = &\ 0.
	\end{aligned}
	$}
	\end{equation*}
	This can then be extended to entire trajectories and the corresponding initial values, i.e.,
	\begin{equation*}
	\resizebox{1\hsize}{!}{$
	\begin{aligned}
	\lim_{k\rightarrow\infty}\lim_{m\rightarrow\infty} \mu \left( \left\{ z_0 \in \mathcal{Z} \mid \| \hat{\varphi}(\eta_0, \cdot, \cdot) - \psi(f(\varphi(y^0, \cdot, \cdot))) \| \geq \epsilon \right\} \right) = 0.
	\end{aligned}
	$}
	\end{equation*}
	Finally, combining this with Theorem~\ref{thm:Koopman_continuous} yields equality for intermediate control values (i.e., $u\in U^p$) for almost all initial conditions $z^0$.\end{proof}
	
\begin{remark} \label{rem:Assumptions}
	Note that in a practical setting, the assumptions $k\rightarrow\infty$ and $m\rightarrow\infty$ cannot be satisfied and moreover, infinitely many data points as well as infinitely large dictionaries are undesirable from a computational point of view. In particular, a low dictionary dimension is crucial to achieve real-time applicability. Nonetheless, the results justify the use of the Koopman operator framework.
\end{remark}

\begin{example}\label{ex:ODE_PBK18}
	Consider the following simple example \cite{BBPK16}:
	\begin{equation}\label{eq:BBPK16}
		\begin{aligned}
			\dot{y}(t) &= G(y(t),u(t)) = \left(\begin{array}{c}
			\mu y_1(t) \\ \lambda (y_2(t) - (y_1(t))^2) + u^{\chi}(t)
			\end{array}\right),\\
			y(0) &= y^0.
		\end{aligned}
	\end{equation}
	Here, we have additionally introduced a factor $\chi$ in order to study the relevance of the linearity assumption for the control input. 
	By restricting ourselves to $n_c$ constant controls, we can transform the system into $n_c$ autonomous systems:
	\begin{equation}\label{eq:BBPK16s}
		\begin{aligned}
			\dot{y}(t) &= G_{u^j}(y(t)) = \left(\begin{array}{c}
			\mu y_1(t) \\ \lambda (y_2(t) - (y_1(t))^2)
			\end{array}\right) + \left(\begin{array}{c}
			0 \\ \left(u^j\right)^{\chi}
			\end{array}\right), \\
			y(0) &= y^0, 
		\end{aligned}
	\end{equation}
	for $j=0,\ldots,n_c-1$.
	\begin{figure}[t]
		\centering
		\parbox[b]{0.24\textwidth}{\centering (a) \\ \includegraphics[width=.24\textwidth]{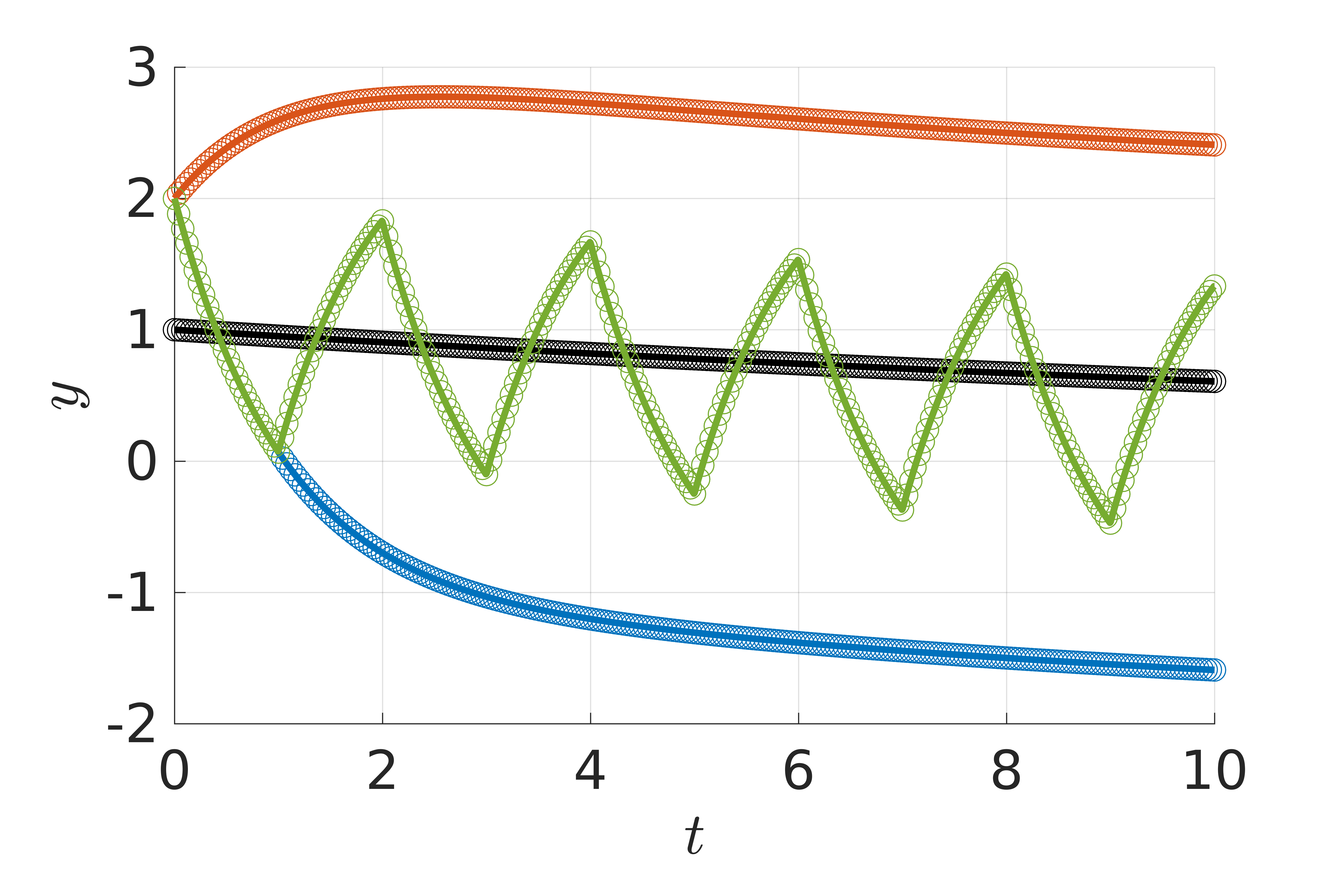}}
		\parbox[b]{0.24\textwidth}{\centering (b) \\ \includegraphics[width=.24\textwidth]{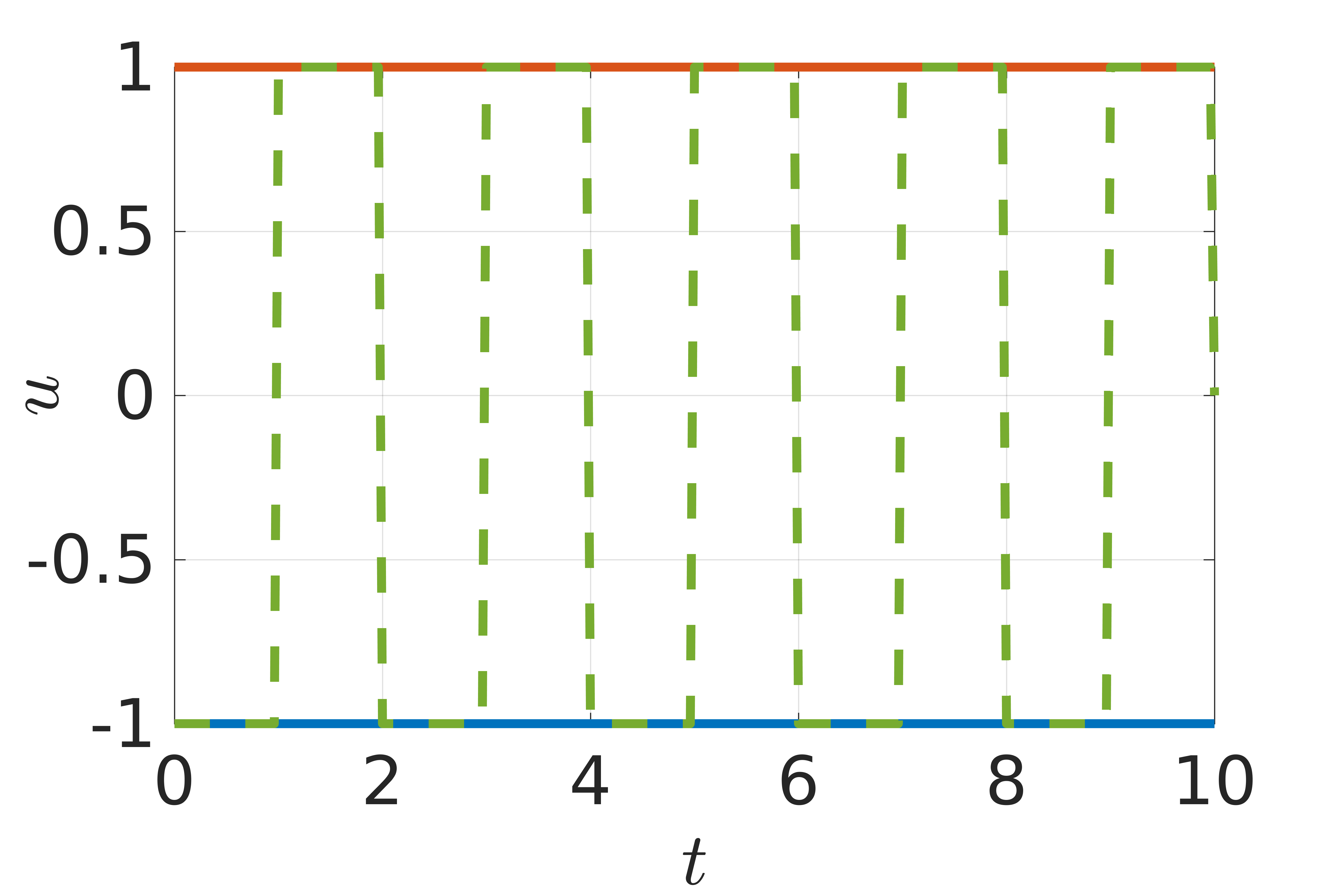}}
		\caption{(a) Trajectories  of $y_2$ of two different autonomous systems ($u^0 = -1$ (blue) and $u^1 = 1$ (orange)) with starting point $ y^0 = (1, 2)^\top $ and linear control input ($\chi = 1$) according to (b). The $y_2$ trajectory of the switched system is shown in green. The lines are the solution of Equation~\eqref{eq:BBPK16s} and the circles the solution of the corresponding K-ROM. The first component $y_1$ is identical in all cases and is shown in black.}
		\label{fig:BBPK_const}
	\end{figure}
	Since this is a finite-dimensional system, we set $f = Id$ (i.e., $z = y$) and observe the full state such that \eqref{eq:KROM_continuous} is a reduced model for $y$. In a first step, we reproduce the dynamics of the two systems with $u^0 = -1$ and $u^1 = 1$ and a switched system with constant switching times (cf.~Figure~\ref{fig:BBPK_const}). For the system~\eqref{eq:BBPK16s}, we can use EDMD to exactly compute the Koopman operator (with $\psi(y) = [y_1, y_2, y_1^2]^\top$, \cite{BBPK16}). Consequently, we observe perfect agreement between the ODE solution and the K-ROM approximation.
	
	In the next step, we compare the solutions of the ODE \eqref{eq:BBPK16} and the K-ROM approximation for intermediate values of $u$, cf.~Figure~\ref{fig:BBPK_continuous}. We see that in accordance with Theorem~\ref{thm:Koopman_continuous}, we again observe a perfect agreement.
	\begin{figure}[h]
		\centering
		\parbox[b]{0.24\textwidth}{\centering (a) \\ \includegraphics[width=.24\textwidth]{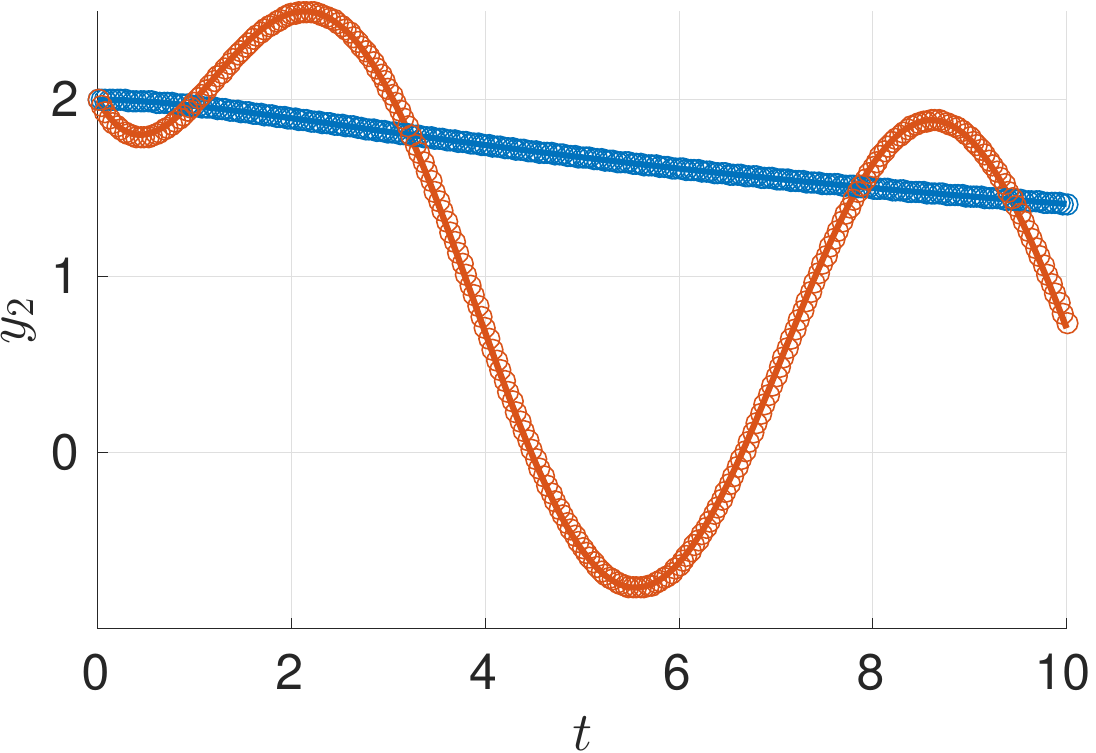}}
		\parbox[b]{0.24\textwidth}{\centering (b) \\ \includegraphics[width=.24\textwidth]{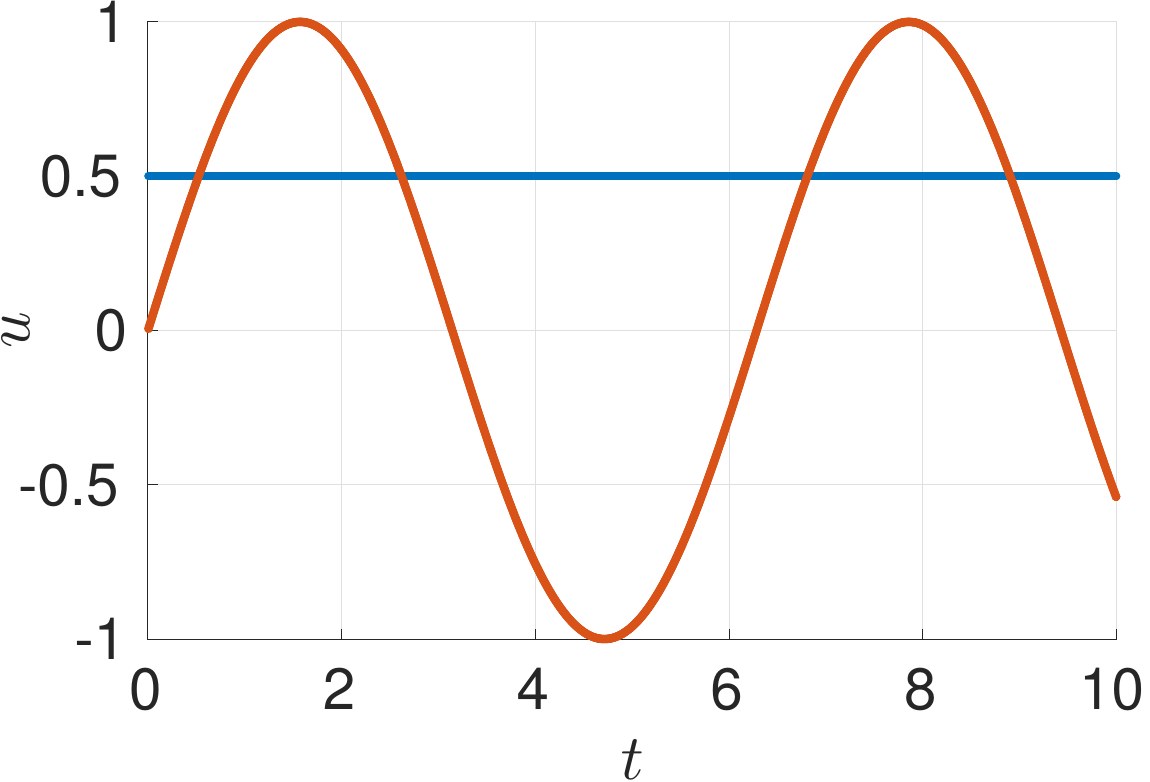}}
		\caption{Similar to Figure~\ref{fig:BBPK_const} but with continuous control inputs, i.e., comparison between \eqref{eq:BBPK16} and the K-ROM approximation.}
		\label{fig:BBPK_continuous}
	\end{figure}
	
	Finally, we study solutions with $\chi > 1$, i.e., nonlinear control inputs. We here construct \eqref{eq:KROM_continuous} from the operators $\Kcal_{0}$ and $\Kcal_{1}$. Due to the choice $u^0 = 0$ and $u^1 = 1$, the exponent $\chi$ does not influence these two solutions and hence, all K-ROM solutions are identical. We observe that the very good agreement from the linear case $\chi = 1$ can not be preserved, cf.~Figure~\ref{fig:BBPK_continuous_nonlinear}, where the $y_2$ trajectories are shown in (a) and the relative error
	\[ \epsilon_{rel}(t) = \frac{ \big| y_2^{\text{\tiny K-ROM}}(t) - y_2^{\text{\tiny ODE}}(t) \big| }{ \big| y_2^{\text{\tiny ODE}}(t) \big| } \]
	is shown in (b). However -- depending on how well the control dependency can be linearized -- we still have acceptable accuracy and a qualitative prediction is possible. Consequently, the K-ROM can be interpreted as being an implicit local linearization with respect to $u$.
	\begin{figure}[h]
		\centering
		\parbox[b]{0.24\textwidth}{\centering (a) \\ \includegraphics[width=.24\textwidth]{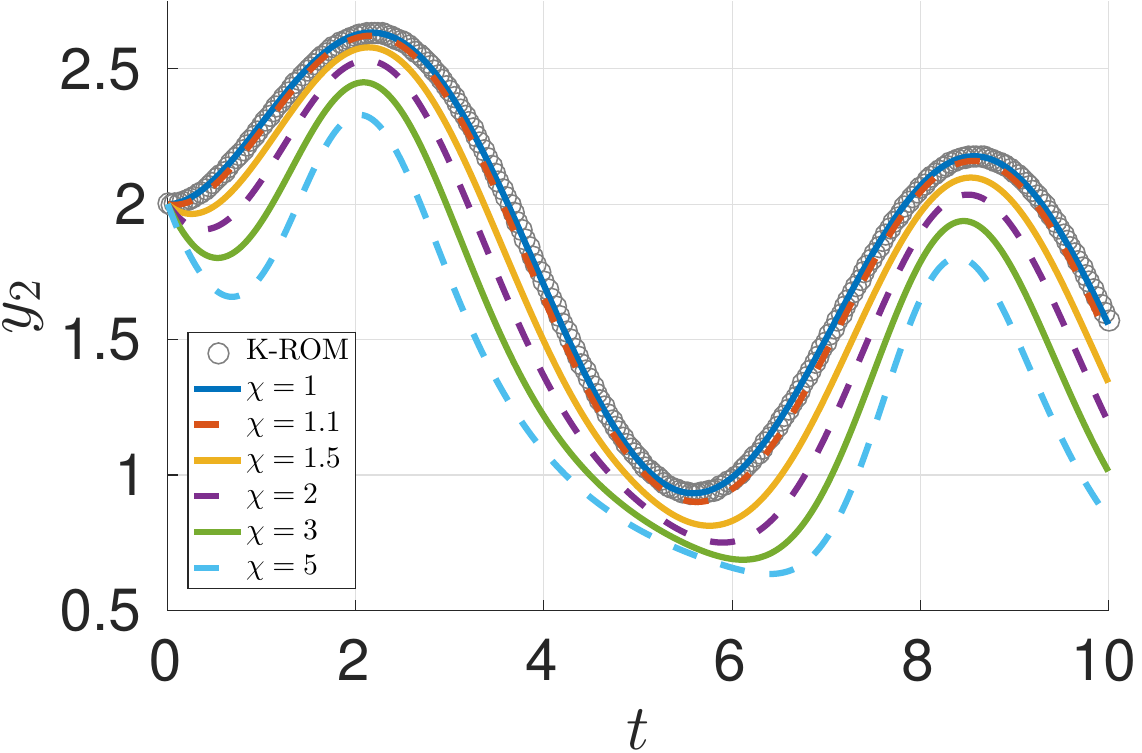}}
		\parbox[b]{0.24\textwidth}{\centering (b) \\ \includegraphics[width=.24\textwidth]{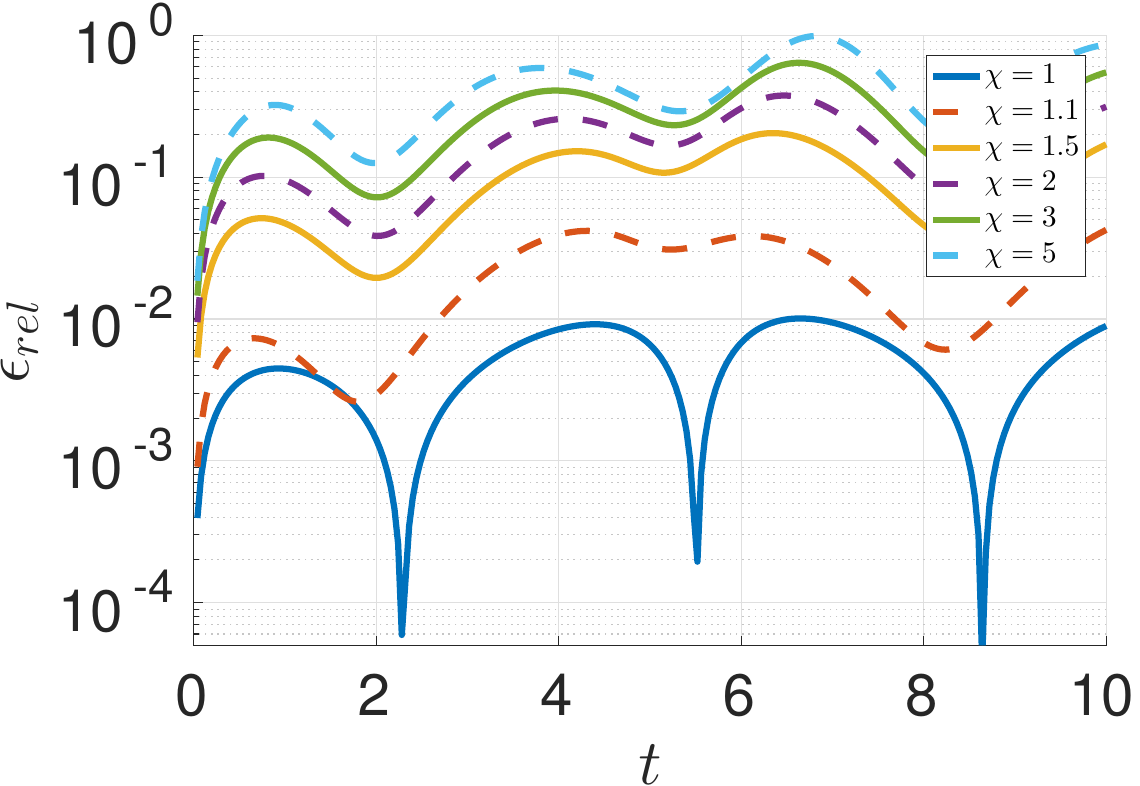}}
		\caption{(a) Comparison between the solutions of the ODE system~\eqref{eq:BBPK16} and \eqref{eq:KROM_continuous} for $u(t) = 0.5 + 0.5 \sin(t)$ and varying values of $\chi$. (b) Relative error between the respective solutions.}
		\label{fig:BBPK_continuous_nonlinear}
	\end{figure}
\end{example}

\section{Control relevant observables}
\label{sec:PDE}
The reduced order model \eqref{eq:KROM_continuous} enables us to approximate a nonlinear control system by a finite-dimensional bilinear system in a very efficient manner. Due to the increased step size of the K-ROM and the linearity, we already observe a significant speed-up. Furthermore, such a bilinear model allows us to use solution methods from bilinear control theory (see \cite{PP08,Ell09} for an introduction), which additionally accelerates the computations. The approach shows its full potential when considering PDEs, which we will show in the following using fluid flow examples with varying complexity. All examples are governed by the 2D incompressible Navier--Stokes equations at moderate Reynolds numbers: 
\begin{align*}
	\dot{y}(x,t) + y(x,t) \cdot \nabla y(x,t) &= \nabla p(x,t) + \frac{1}{Re} \Delta y(x,t), \\
	\nabla \cdot y(x,t) &= 0, \\
	y(x,0) &= y^0(x),
\end{align*}
where $y$ is the flow velocity and $p$ is the pressure, and with appropriately chosen boundary conditions. The problem is discretized using a finite volume method and the PISO scheme for time integration \cite{FP02}. All calculations are performed using the open source solver OpenFOAM \cite{JJT07}.

\begin{figure}[b]
	\centering
	\parbox[b]{0.24\textwidth}{\centering (a) \\ \includegraphics[width=.24\textwidth]{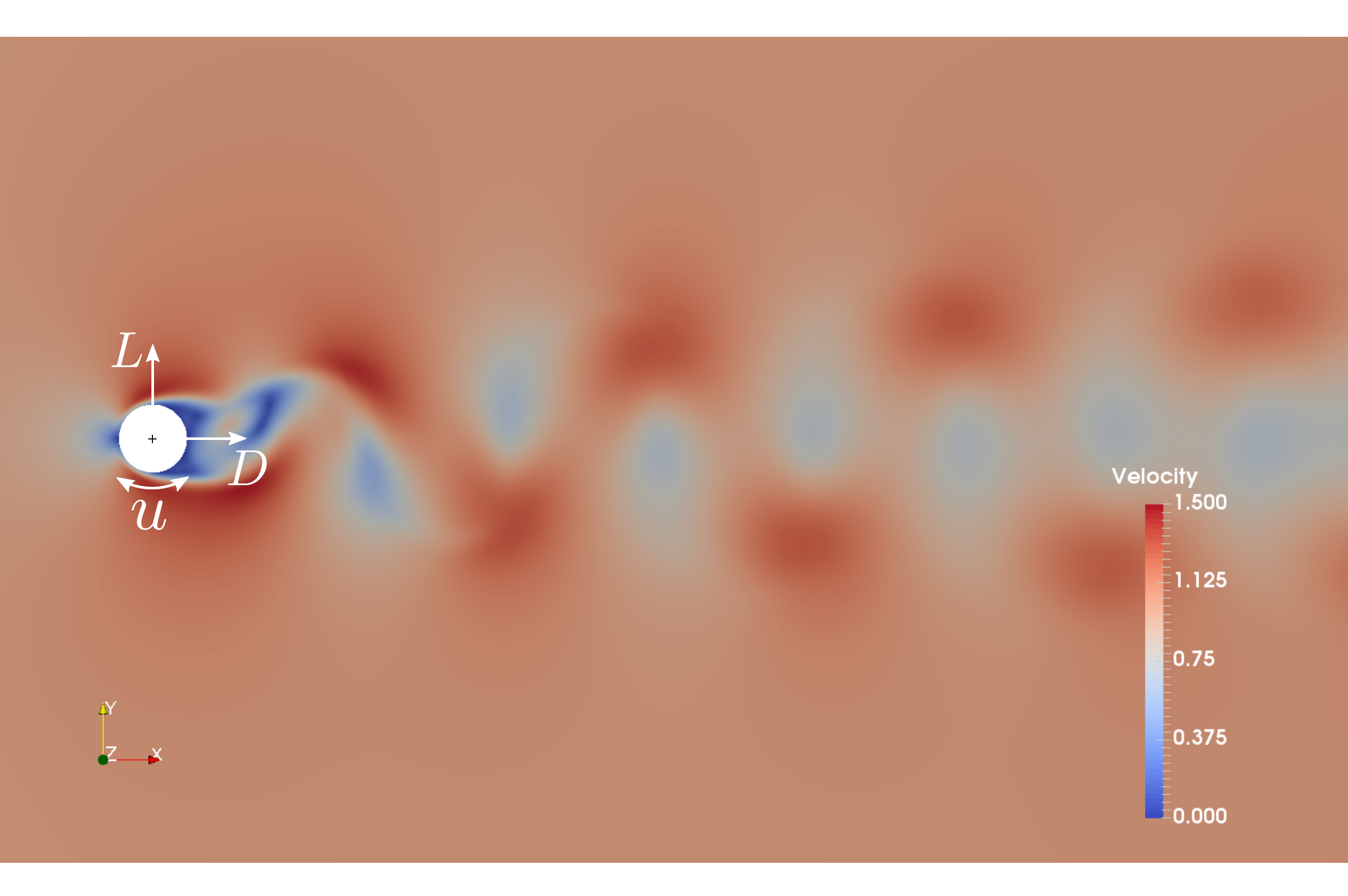}}
	\parbox[b]{0.24\textwidth}{\centering (b) \\ \includegraphics[width=.24\textwidth]{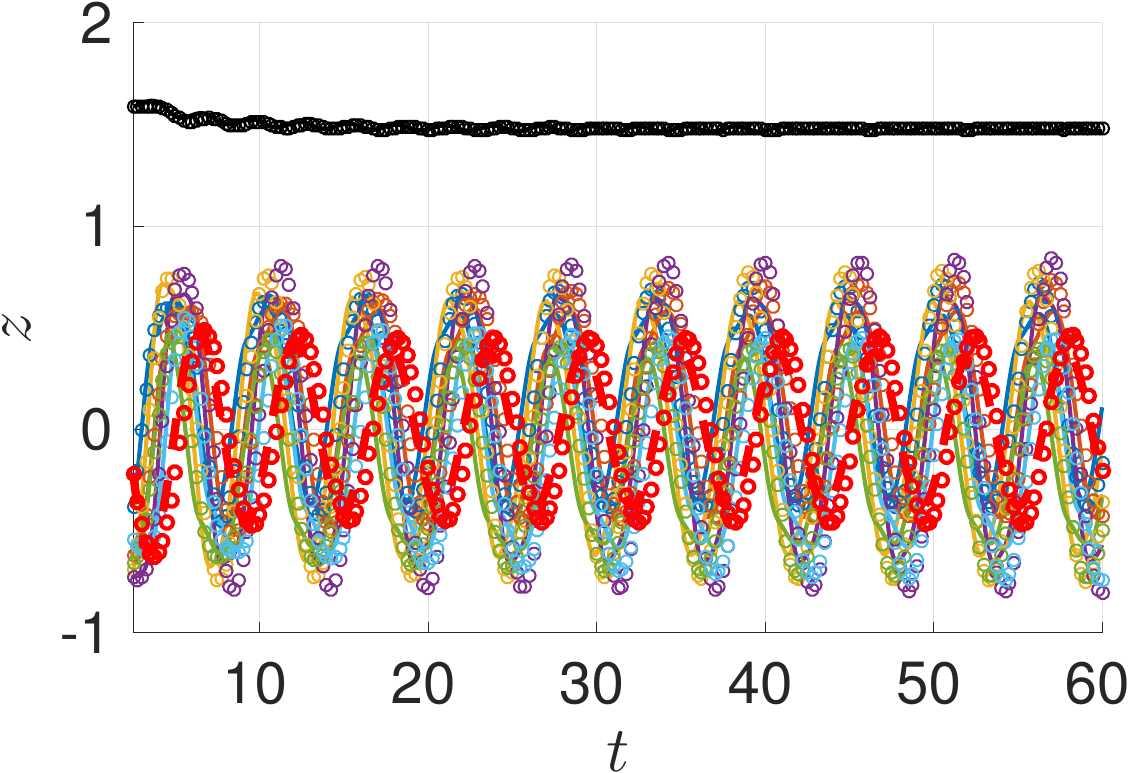}} \\[1ex]
	\parbox[b]{0.24\textwidth}{\centering (c) \\ \includegraphics[width=.24\textwidth]{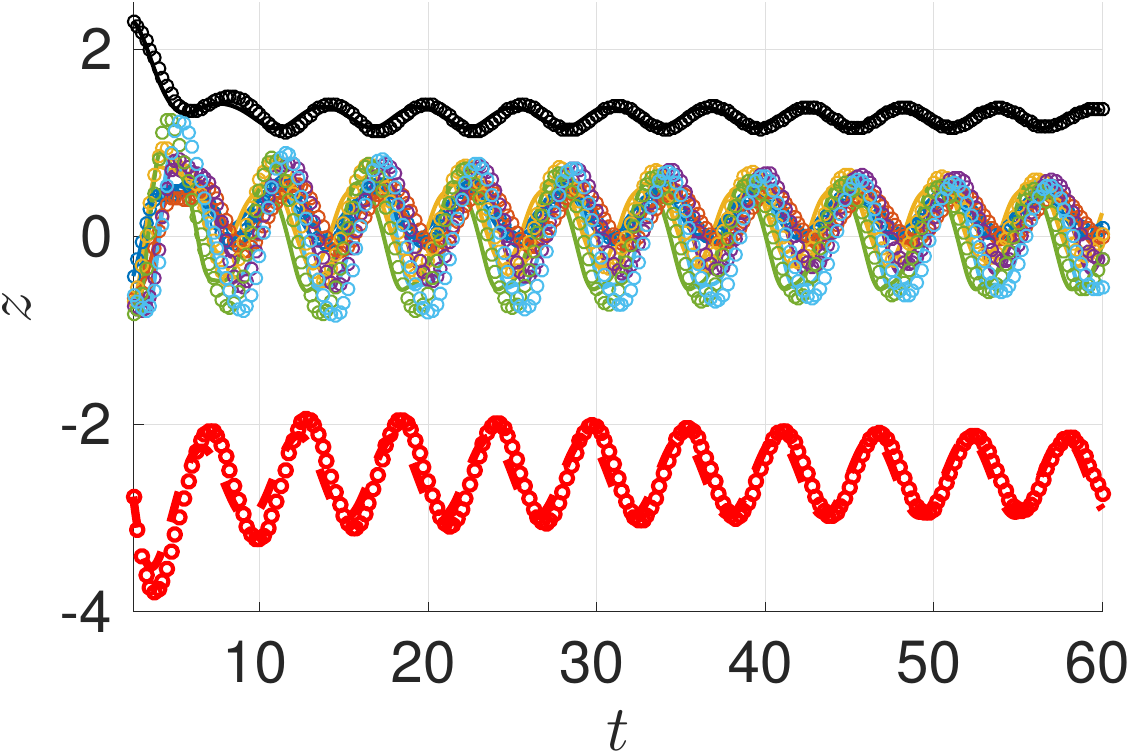}}
	\parbox[b]{0.24\textwidth}{\centering (d) \\ \includegraphics[width=.24\textwidth]{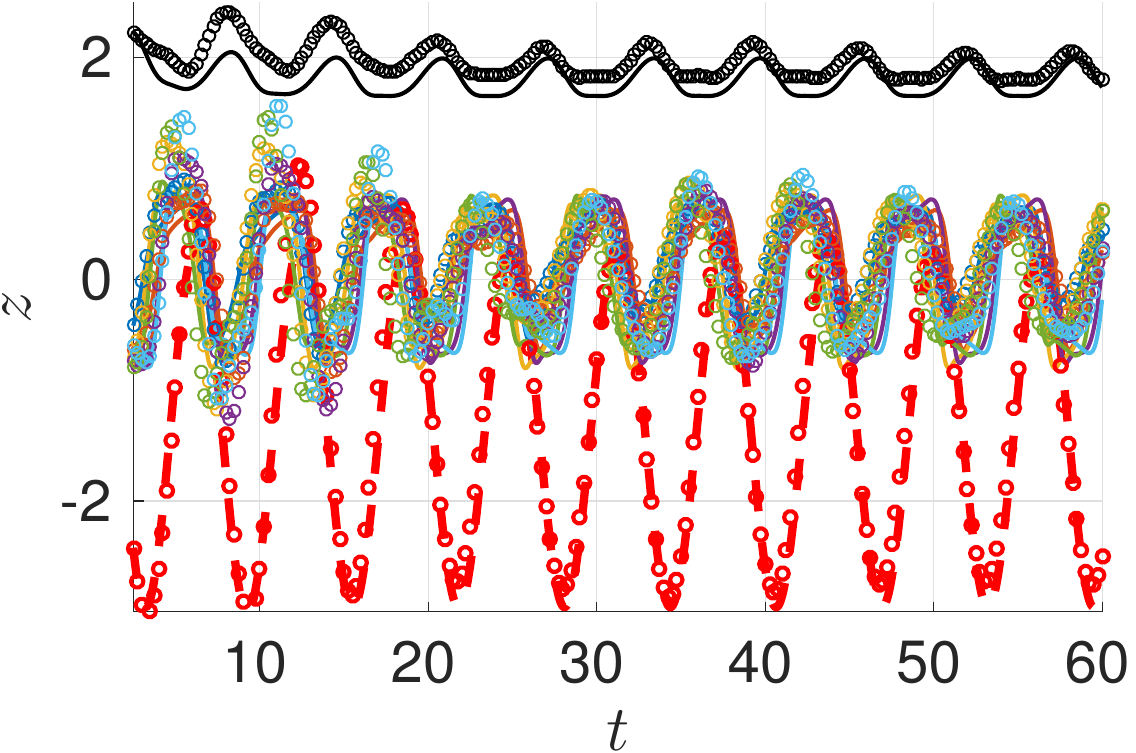}}
	\caption{(a) Sketch of the problem setting. The system is controlled by rotating the cylinder. (b) and (c) Observation of the PDE solution and K-ROM approximation for $u^a = 0$ and $u^b = 2$, respectively. The lift is shown by the dashed red line, the drag is shown in black and the remaining lines are the vertical velocities at the six positions shown in (a). (d) Observation of the PDE solution and K-ROM approximation for $u(t) = 1 + sin(t)$.}
	\label{fig:vonKarman}
\end{figure}

As a first example, we consider the widely studied flow around a cylinder at a Reynolds number of $Re = 100$ (see Figure~\ref{fig:vonKarman}~(a) for the problem setup). The system is controlled via rotation of the cylinder, i.e., $u(t)$ is the angular velocity. The uncontrolled system possesses a periodic solution, the well-known \emph{von K\'{a}rm\'{a}n vortex street}.

In order to obtain a low-dimensional reduced model and achieve a high acceleration, we construct the K-ROM only for the control relevant quantities, i.e., the lift  $L$ and drag force $D$ of the cylinder acting in vertical and horizontal direction, respectively. Additionally, we observe the vertical velocity at six different positions $(x_1,\ldots,x_6)$ in the cylinder wake (see Figure~\ref{fig:vonKarman}~(a)):
\begin{align*}
	z_i &= f((y(\cdot,t_i), p(\cdot,t_i)) \\&= \left(L(t_i), D(t_i), y_2(x_1,t_i), \ldots, y_2(x_6,t_i)\right)^{\top}.
\end{align*}

For the approximation via EDMD, we use as basis functions all monomials up to order 2, which results in a K-ROM of dimension 45. We collect data from one long term simulation with random control inputs $u_i\in\hat{U}=\{0,2\}$, which is then split into two parts corresponding to the two autonomous systems. The time step in the time series is $0.25$ in comparison to the time step $0.01$ in the finite volume scheme. Figures~\ref{fig:vonKarman}~(b) and \ref{fig:vonKarman}~(c) show a comparison between the PDE and the K-ROM solution for constant control inputs $u^0 = 0$ and $u^1 = 2$. We see that in both cases, the solutions agree remarkably well, considering that the finite volume discretization for this case consists of $22,000$ cells (i.e., $66,000$ unknowns per time step) and the K-ROM is a 45-dimensional linear model. This results in a speed-up of approximately $75,000$ (OpenFOAM vs.~MATLAB). Figure~\ref{fig:vonKarman}~(d) shows a comparison of the solutions for a sinusoidal control and we see that the agreement is still satisfactory although the system does not depend linearly on the control such that the assumptions of Theorem~\ref{thm:Koopman_continuous} are violated. In the next section, we will use the K-ROM approach to solve a more difficult flow control problem.

\section{Model predictive control}
\label{sec:MPC}
We have seen in the previous section that we can use \eqref{eq:KROM_continuous} to approximate a nonlinear, infinite-dimensional control system by a finite-dimensional bilinear system in a very efficient manner. However, since the accuracy cannot always be expected to be this good for finite amounts of data and finite-dimensional dictionaries $\psi$ -- in particular for more complex dynamics -- we now embed the above approach in an MPC framework such that we only require high accuracy on short time intervals. To this end, we formulate the K-ROM approximation of the original closed-loop problem \eqref{eq:MPC}:
\begin{align}
	\min_{{u} \in U^p} \sum_{i=s}^{s+p-1} &\hat{L}(\eta_i) + \alpha \|u_i\|_2 + \beta \| u_i - u_{i-1} \|_2 \label{eq:MPC_Koopman}\tag{K-MPC}\\
	\text{s.t.}\quad \eta_{i+1} &= A \eta_{i} + \sum_{j=1}^{n_c-1} B_j \eta_{i} \frac{u_{j,i} - u^0}{u^j-u^0}, \notag \\
	\eta_s &= \psi(f(y_s)), \notag
\end{align}
The K-ROM based MPC method is summarized in Algorithm~\ref{alg:MPC_Koopman}.
\begin{algorithm}[t]
	\caption{(K-ROM-based MPC)}
	\label{alg:MPC_Koopman}
	\begin{algorithmic}[1]
		\Require EDMD approximations of $n_c$ Koopman operators; prediction horizon length $p \in \N$.
		\For{$i=0, 1, 2, \ldots$}.
		\State \parbox[t]{\dimexpr\linewidth-\algorithmicindent}{Observe current state: $\eta_i = \psi(z_i) = \psi(f(y_i))$. \strut}
		\State \parbox[t]{\dimexpr\linewidth-\algorithmicindent}{Solve Problem~\eqref{eq:MPC_Koopman} with initial condition $\eta_{i}$ on the prediction horizon of length $p$. \strut}
		\State \parbox[t]{\dimexpr\linewidth-\algorithmicindent}{Apply the first entry, i.e., $u^*_1$, to the system. \strut}
		\EndFor
	\end{algorithmic}
\end{algorithm}
Equality of the objective functions follows from Theorem~\ref{thm:Koopman_continuous} and Corollary~\ref{cor:EqualDynamics}.
\begin{theorem}\label{thm:EqualMPC} 
	Consider Problem \eqref{eq:MPC} and assume that we have convergence of EDMD towards the Koopman operator according to \cite{KM18a}. 
	Furthermore, assume $L(y(t)) = \hat{L}(\eta_{i})$ for all $t\in[t_0,t_e]$ and the corresponding $i = (t-t_0)/h$.
	Then, as the basis size $k$ and number of sampled data points $m$ tend to infinity,
	the objective function value of Problem \eqref{eq:MPC_Koopman} converges in measure (with respect to the initial condition $z_0 = f(y^0)$) to that of Problem \eqref{eq:MPC} for every $u \in U^p$, $p < \infty$. That is, for all $\epsilon > 0$ we obtain
	\vspace{-.25em}
	\[
	\lim_{k\rightarrow\infty} \lim_{m\rightarrow\infty} \mu \Big( \Big\{ z_0 \in \mathcal{Z} \Big | \sum_{i=0}^{p} | \hat{L}(\eta_i) - L(y_i) | \geq \epsilon \Big\} \Big) = 0.
	\]
\end{theorem}
Since the focus of this article is on the development of the K-ROM, we will simply solve Problem~\eqref{eq:MPC_Koopman} using MATLAB's internal SQP solver. Nevertheless, an additional increase in efficiency can likely be obtained using methods known from bilinear systems theory \cite{PP08}.

\begin{figure}[h]
	\centering
	\includegraphics[width=.48\textwidth]{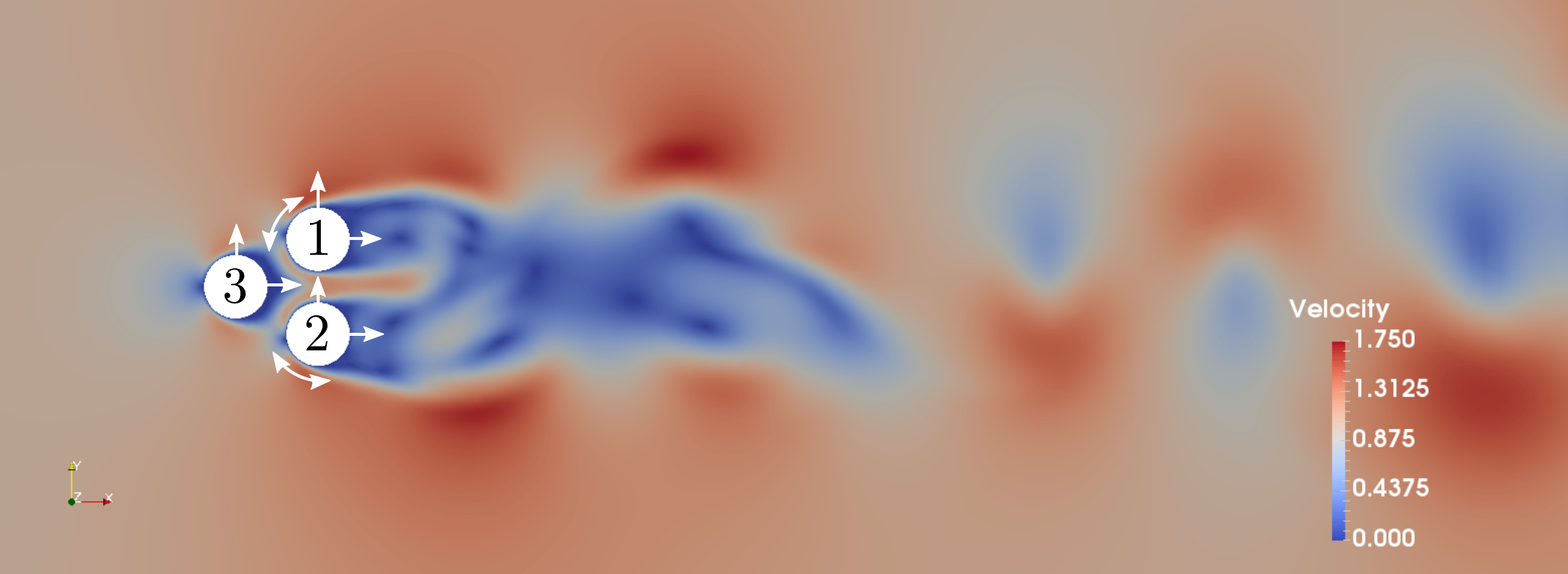}
	\caption{Setup of the fluidic pinball according to \cite{DPMN18}}
	\label{fig:fluidicPinball}
\end{figure}
To demonstrate the effectiveness of the K-ROM MPC approach, we again use a flow control problem, but of higher complexity. We now consider the \emph{fluidic pinball} \cite{DPMN18}, which is a configuration with three cylinders placed on the edges of an equidistant triangle with edge length $1.5 R$, where $R$ is the cylinder radius, cf.~Fig.~\ref{fig:fluidicPinball}. In order to entirely avoid point-wise evaluations of the flow field, we only observe the forces acting on the cylinders. To increase the number of observables, we use delay coordinates, which is a very common approach (see, for instance, \cite{BBP+17}):
\begin{equation*}
\resizebox{1\hsize}{!}{$
	\begin{aligned}
	z_i = \left(L^1_i, L^2_i, L^3_i, D^1_i, D^2_i, D^3_i, L^1_{i-1}, L^2_{i-1}, L^3_{i-1}, D^1_{i-1}, D^2_{i-1}, D^3_{i-1}\right)^\top
	\end{aligned}
$}
\end{equation*}
The goal is to control the lift forces of all three cylinders by rotating the cylinders two and three. Hence, we simply have to track the corresponding entries of $z$ in the MPC problem:
\begin{equation*}
	\resizebox{1\hsize}{!}{$
	\begin{aligned}
	\min_{u \in [0,1]^{2\times p}} \sum_{i=s}^{s+p-1} &\left(\sum_{j=1}^{3}\left(z_{j,i} - z_{j}^{\mathsf{ref}}\right)^2\right)  + \alpha \|u_i\|_2 + \beta \| u_i - u_{i-1} \|_2\\
	\text{s.t.}\quad \eta_{i+1} = &A \eta_{i} + \sum_{j=1}^{2} B_j \eta_{i} \hat{u}_{j,i} \quad \text{for}~i = s,\ldots,s + p - 1, \\
	\eta_s = &\psi(f(y_s, p_s)).
	\end{aligned}
	$}
\end{equation*}
We want to allow control inputs between $-2$ and $2$ in each direction and since the control does not enter linearly into the system, we create multiple K-ROMs which are \emph{localized} in the control domain. This means that we approximate nine Koopman operators $\Kcal_{\mathcal{I}}$ with 
\begin{equation*}
\resizebox{1\hsize}{!}{$
\mathcal{I} \in \{ [-2,-2], [-2,0], [-2,2], [0,-2],[0,0], [0,2], [2,-2], [2,0], [2,2] \}
$}
\end{equation*}
and construct six reduced models of the form \eqref{eq:KROM_continuous}, each of which is valid on the respective part of the control domain. Similar ideas of localized reduced order models are often used in the \emph{reduced basis} community (cf., e.g., \cite{AHKO12,BDPV18}).

We now study the performance of the behavior of the K-ROM based controller for different Reynolds numbers. In \cite{DPMN18} it was argued that the fluidic pinball without control possesses a periodic solution for $Re<90$, then a quasi-periodic solution, and that the system behaves chaotically for $Re \geq 120$. Consequently, it becomes more and more challenging to construct accurate surrogate models with increasing Reynolds number.

\begin{figure}[h]
	\centering
	\parbox[b]{0.24\textwidth}{\centering (a) \\ \includegraphics[width=.24\textwidth]{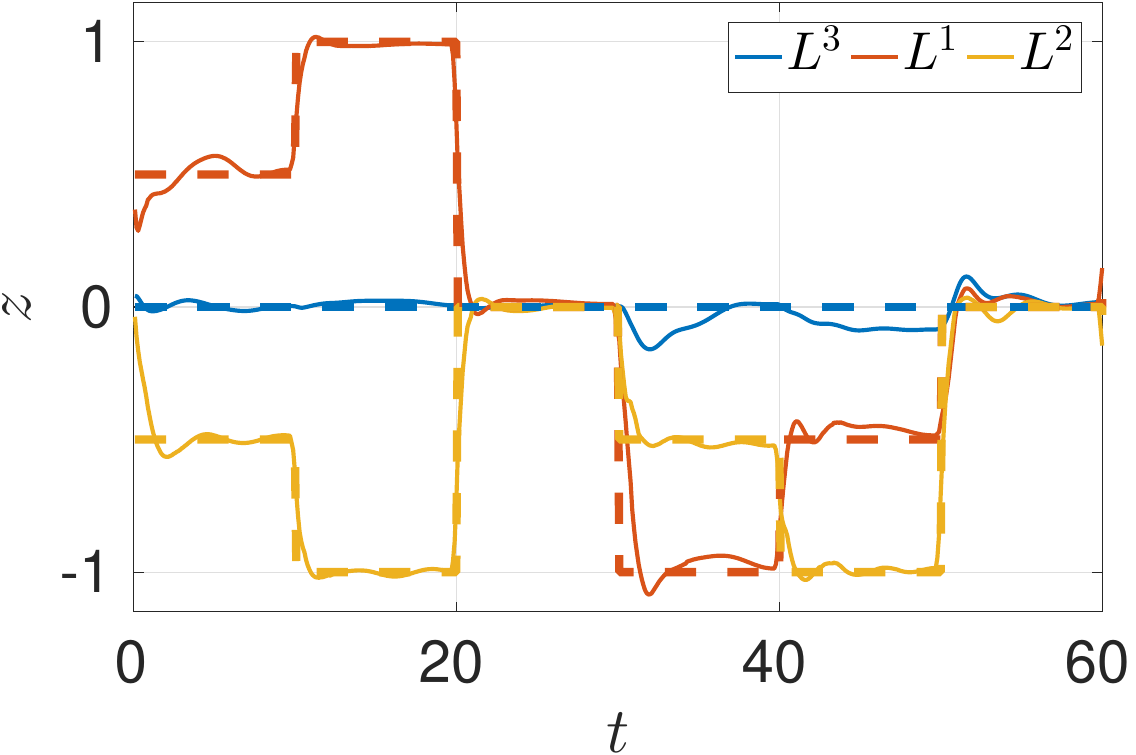}}
	\parbox[b]{0.24\textwidth}{\centering (b) \\ \includegraphics[width=.24\textwidth]{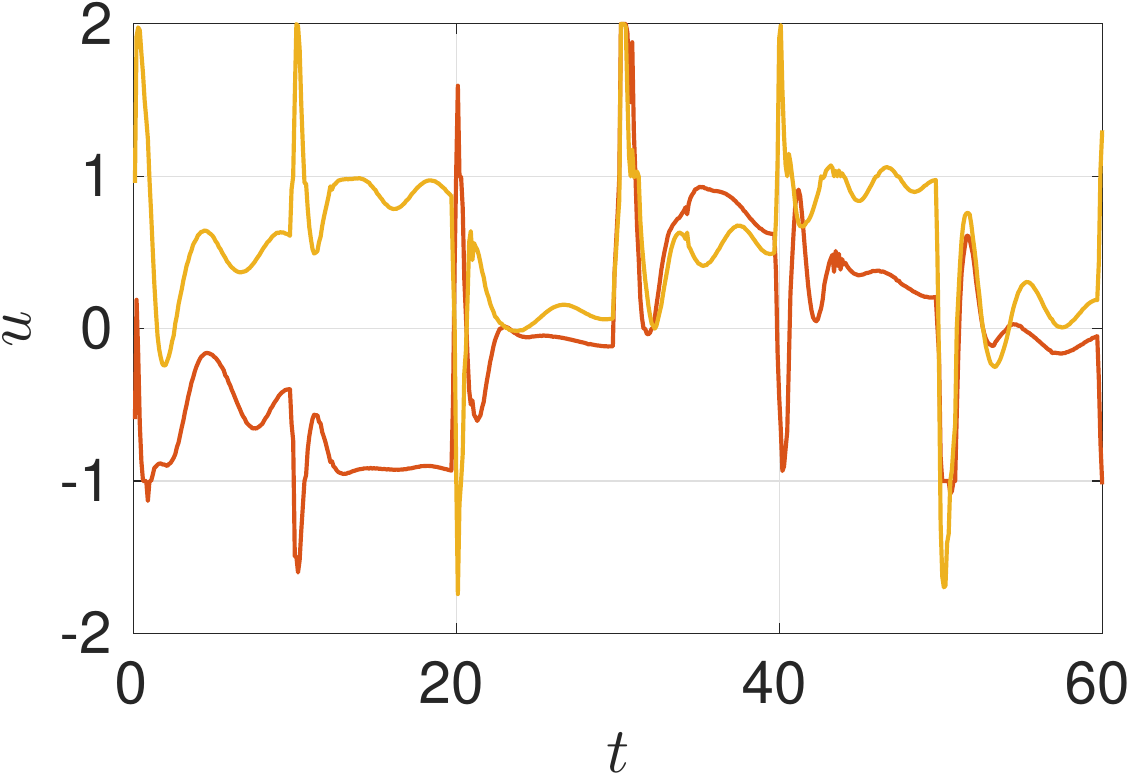}}
	\caption{(a) Lift and corresponding reference trajectories for the three lift forces at $Re = 100$, computed with Algorithm~\ref{alg:MPC_Koopman} with $p=5$. (b) The corresponding control inputs for the cylinders 1 and 2.}
	\label{fig:fluidicPinball_Re100}
\end{figure}

As the first case, we set $Re=100$, i.e., we have quasi-periodic dynamics in the uncontrolled system. In this case, a simple DMD approximation ($\psi(z) = z$) is sufficiently accurate, and we obtain a twelve-dimensional bilinear surrogate model with a time step of $0.1$. In comparison to solving the full system (where we use adaptive time stepping such that the CFL number is below 0.5), a speed-up factor of more than six orders of magnitude is obtained, and real-time applicability is achieved. The resulting system behavior for a piecewise constant reference state is shown in Fig.~\ref{fig:fluidicPinball_Re100} and we observe excellent control performance.

\begin{figure}[h]
	\centering
	\parbox[b]{0.24\textwidth}{\centering (a) \\ \includegraphics[width=.24\textwidth]{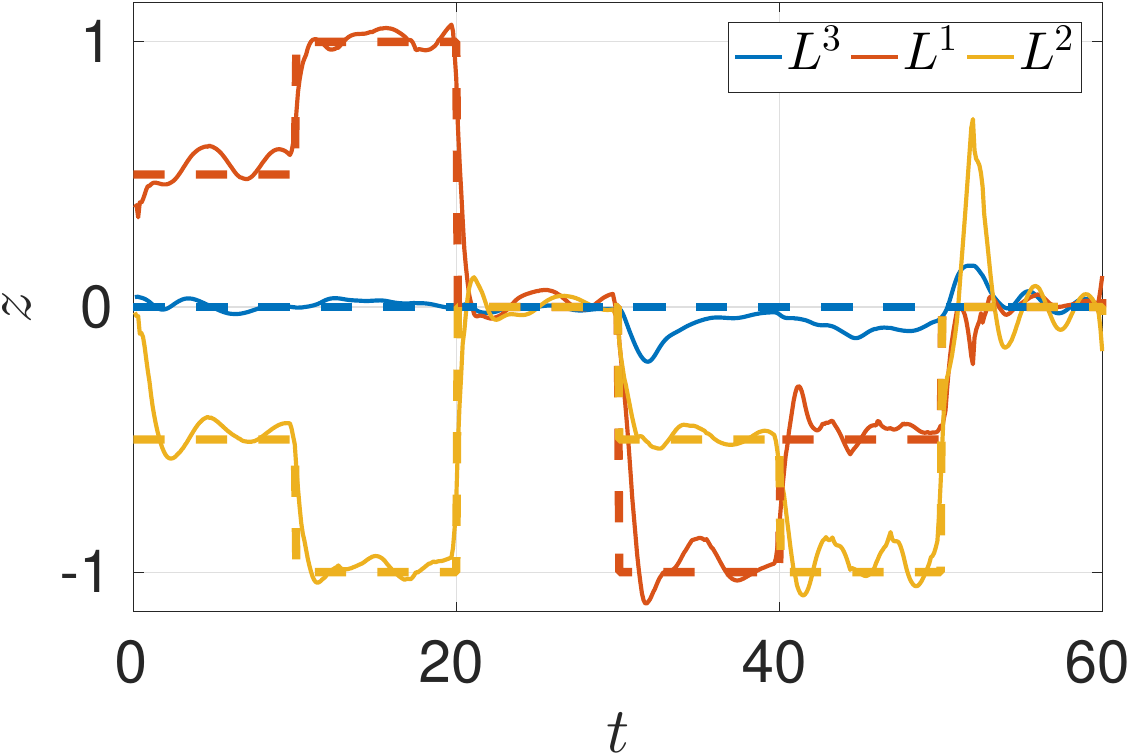}}
	\parbox[b]{0.24\textwidth}{\centering (b) \\ \includegraphics[width=.24\textwidth]{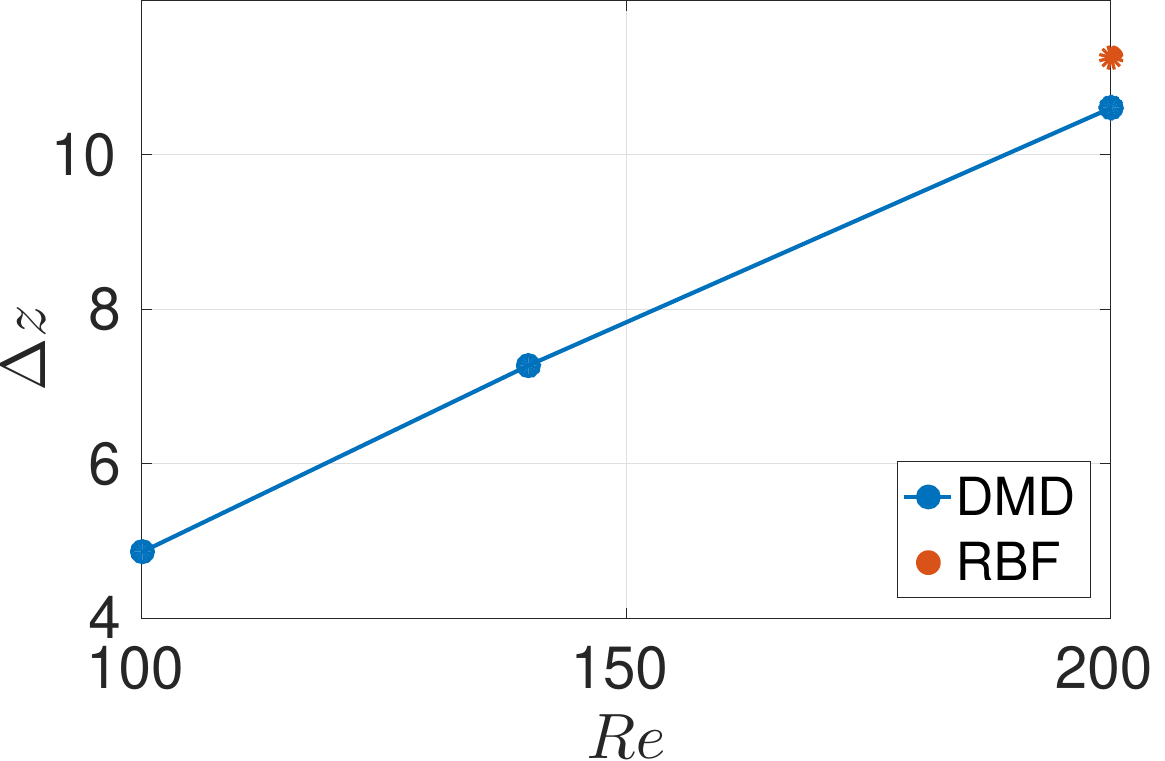}}
	\caption{(a) Similar to Fig.~\ref{fig:fluidicPinball_Re100} (a) but with $Re = 140$. (b) The control error $\Delta z$ with increasing Reynolds number.}
	\label{fig:fluidicPinball_Re140}
\end{figure}

As the second case, we set $Re=140$ such that the system behaves mildly chaotically. We see in Fig.~\ref{fig:fluidicPinball_Re140}~(a) that the fluctuations around the desired state are increased. This is due to two facts. First, the system is chaotic and hence much more difficult to control. Furthermore, the accuracy of the K-ROM is lower than for the quasi-periodic case such that the MPC problems \eqref{eq:MPC} and \eqref{eq:MPC_Koopman} do not necessarily possess the same solutions any longer. Nevertheless, the control task is still performed quite satisfactorily, considering that we have replaced a nonlinear PDE with $\approx 50,000$ cells and $\approx 150,000$ degrees of freedom by a twelve-dimensional linear system.

\begin{figure}[h]
	\centering
	\parbox[b]{0.24\textwidth}{\centering (a) \\ \includegraphics[width=.24\textwidth]{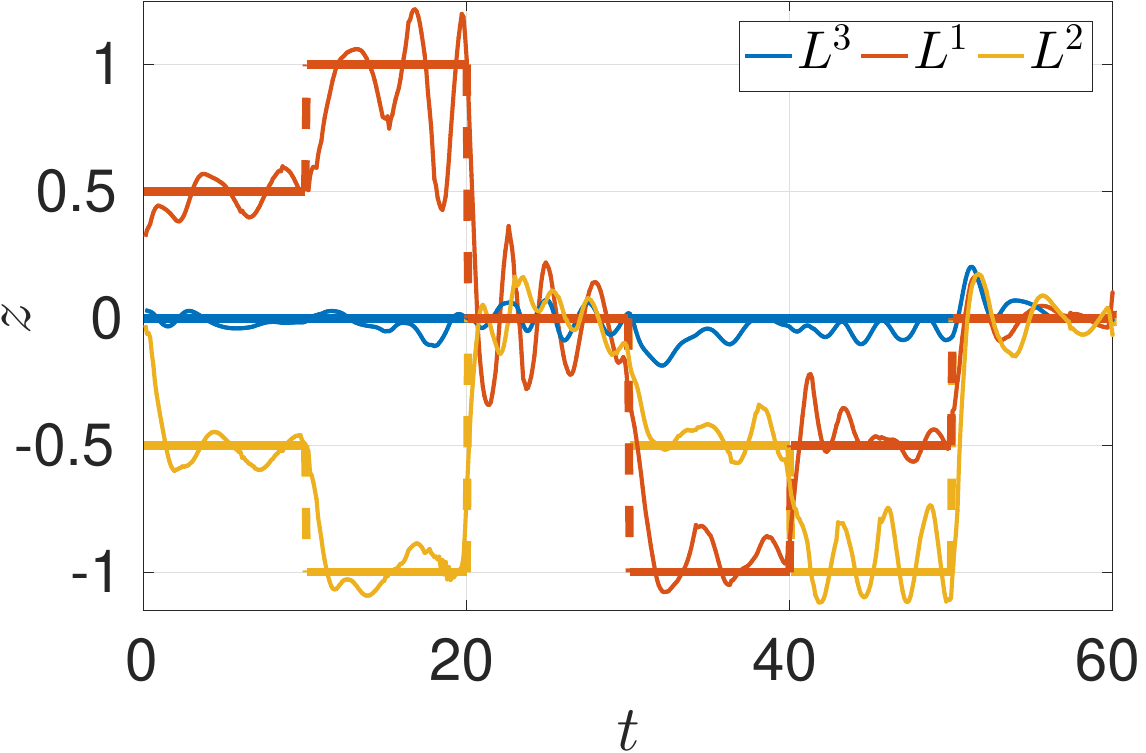}}
	\parbox[b]{0.24\textwidth}{\centering (b) \\ \includegraphics[width=.24\textwidth]{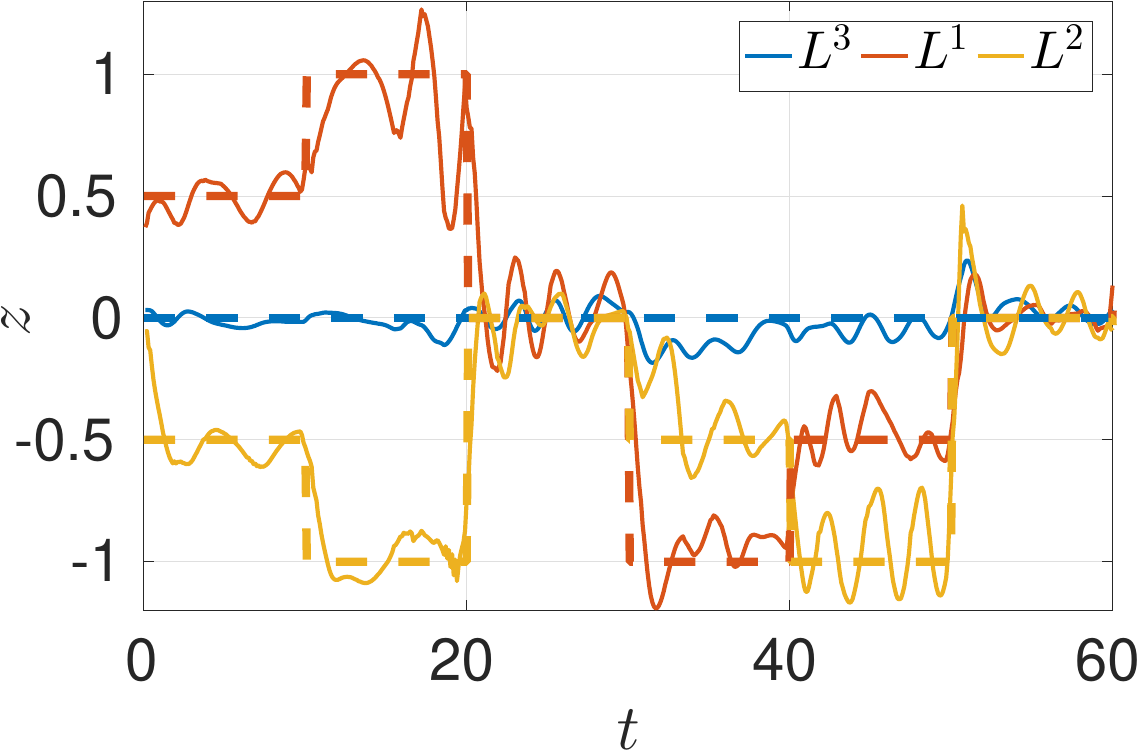}}
	\caption{Similar to Fig.~\ref{fig:fluidicPinball_Re100} (a) but with $Re = 200$. (a) Approximation via DMD as for the tow previous cases. (b) Approximation via 1000 randomly distributed radial basis functions.}
	\label{fig:fluidicPinball_Re200}
\end{figure}

When further increasing the Reynolds number well into the chaotic regime, we observe still stronger oscillations around the desired state, see also Fig.~\ref{fig:fluidicPinball_Re140}~(b), where the tracking error $\Delta z = \int_{0}^{60}  \left \| z(t) - z^{\mathsf{ref}(t)}\right \|^2 \, dt$ is shown. This is due to the two facts mentioned above, complexity of the dynamics and prediction quality. In order to study the dependence of the prediction quality, we compare the approximation via DMD to one where we additionally place $1000$ radial basis functions in the 12-dimensional observation space as a Halton set (i.e., quasi-randomly). The two runs are compared in Fig.~\ref{fig:fluidicPinball_Re200} and also in Fig.~\ref{fig:fluidicPinball_Re140}~(b). We see that the much higher dimension of $\psi$ has no positive impact on the control performance, which indicates that the inaccuracy is mainly due to the chaotic dynamics.

\section{Conclusion}
\label{sec:Conclusion}

We have presented a new approach for Koopman operator based reduced order models which can be used for real-time control of nonlinear PDEs. Multiple Koopman operators are approximated via EDMD at different constant control inputs, and intermediate control values can be approximated by linear interpolation between these operators which yields a bilinear control system. For control affine systems, we obtain convergence in measure of the reduced towards the full objective function. However, as examples show, the approach leads to very good results even in situations where this assumption does not hold. In this case, the control approach can be interpreted as an implicit local linearization.

Due to the larger step sizes and the linearity of the K-ROM, the reduced model can be solved significantly faster, in the case of the 2D Navier--Stokes equations by five to six orders of magnitude. An additional benefit is that since the model is bilinear, we can use efficient solution methods for the reduced control problem. Due to the restriction to several autonomous systems, the training data requirements are very low.

One further direction of research is to develop stronger statements about the error for the K-ROM approach, e.g., concerning the basis size or the required data. To this end, it might also be interesting to study the influence of the assumptions on the dynamical system in order to obtain convergence EDMD towards the Koopman operator. From a control theoretic perspective, it would be very interesting to investigate whether the notion of controllability can be carried over to nonlinear systems. Moreover, feedback controllers for bilinear systems could help to further improve the efficiency over MPC. In terms of numerical efficiency, automated methods for choosing appropriate basis functions for the system dynamics (e.g., via dictionary learning \cite{LDBK17}) could help to further improve the range of applicability. 

\bibliographystyle{unsrt}
\bibliography{Bibliography}

\begin{thebibliography}{10}

\bibitem{PK19}
S.~Peitz and S.~Klus.
\newblock {Koopman operator-based model reduction for switched-system control
  of PDEs}.
\newblock {\em Automatica (to appear)}, 55(8), 2019.

\bibitem{BN15}
S.~L. Brunton and B.~R. Noack.
\newblock {Closed-Loop Turbulence Control: Progress and Challenges}.
\newblock {\em Applied Mechanics Reviews}, 67(5):1--48, 2015.

\bibitem{GP17}
L.~Gr{\"{u}}ne and J.~Pannek.
\newblock {\em {Nonlinear Model Predictive Control}}.
\newblock Springer International Publishing, 2 edition, 2017.

\bibitem{DBN17}
T.~Duriez, S.~L. Brunton, and B.~R. Noack.
\newblock {\em {Machine Learning Control – Taming Nonlinear Dynamics and
  Turbulence}}.
\newblock Springer, 2017.

\bibitem{LMQR14}
T.~Lassila, A.~Manzoni, A.~Quarteroni, and G.~Rozza.
\newblock {Model order reduction in fluid dynamics : challenges and
  perspectives}.
\newblock In Alfio Quarteroni and Gianluigi Rozza, editors, {\em Reduced Order
  Methods for Modeling and Computational Reduction}, pages 235--273. Springer,
  Cham, 2014.

\bibitem{BGW15}
P.~Benner, S.~Gugercin, and K.~Willcox.
\newblock {A Survey of Projection-Based Model Reduction Methods for Parametric
  Dynamical Systems}.
\newblock {\em SIAM Review}, 57(4):483--531, 2015.

\bibitem{Sir87}
L.~Sirovich.
\newblock {Turbulence and the dynamics of coherent structures part I: coherent
  structures}.
\newblock {\em Quarterly of Applied Mathematics}, XLV(3):561--571, 1987.

\bibitem{KV99}
K.~Kunisch and S.~Volkwein.
\newblock {Control of the Burgers Equation by a Reduced-Order Approach Using
  Proper Orthogonal Decomposition}.
\newblock {\em Journal of Optimization Theory and Applications},
  102(2):345--371, 1999.

\bibitem{Row05}
C.~W. Rowley.
\newblock {Model Reduction for Fluids, Using Balanced Proper Orthogonal
  Decomposition}.
\newblock {\em International Journal of Bifurcation and Chaos},
  15(3):997--1013, 2005.

\bibitem{HV05}
M.~Hinze and S.~Volkwein.
\newblock {Proper Orthogonal Decomposition Surrogate Models for Nonlinear
  Dynamical Systems: Error Estimates and Suboptimal Control}.
\newblock In P.~Benner, D.~C. Sorensen, and V.~Mehrmann, editors, {\em
  Reduction of Large-Scale Systems}, volume~45, pages 261--306. Springer Berlin
  Heidelberg, 2005.

\bibitem{BDPV18}
D.~Beermann, M.~Dellnitz, S.~Peitz, and S.~Volkwein.
\newblock {Set-Oriented Multiobjective Optimal Control of PDEs using Proper
  Orthogonal Decomposition}.
\newblock In {\em Reduced-Order Modeling (ROM) for Simulation and
  Optimization}, pages 47--72. Springer, 2018.

\bibitem{Koo31}
B.~O. Koopman.
\newblock {Hamiltonian Systems and Transformations in Hilbert Space}.
\newblock {\em Proceedings of the National Academy of Sciences},
  17(5):315--318, 1931.

\bibitem{Mez05}
I.~Mezi{\'{c}}.
\newblock {Spectral Properties of Dynamical Systems, Model Reduction and
  Decompositions}.
\newblock {\em Nonlinear Dynamics}, 41:309--325, 2005.

\bibitem{BMM12}
M.~Budi{\v{s}}i{\'{c}}, R.~Mohr, and I.~Mezi{\'{c}}.
\newblock {Applied Koopmanism}.
\newblock {\em Chaos}, 22, 2012.

\bibitem{Sch10}
P.~J. Schmid.
\newblock Dynamic mode decomposition of numerical and experimental data.
\newblock {\em Journal of Fluid Mechanics}, 656:5--28, 2010.

\bibitem{RMB+09}
C.~W. Rowley, I.~Mezi{\'{c}}, S.~Bagheri, P.~Schlatter, and D.~S. Henningson.
\newblock {Spectral analysis of nonlinear flows}.
\newblock {\em Journal of Fluid Mechanics}, 641:115--127, 2009.

\bibitem{TRL+14}
J.~H. Tu, C.~W. Rowley, D.~M. Luchtenburg, S.~L. Brunton, and J.~N. Kutz.
\newblock {On Dynamic Mode Decomposition: Theory and Applications}.
\newblock {\em Journal of Computational Dynamics}, 1(2):391--421, 2014.

\bibitem{KGPS18}
S.~Klus, P.~Gel{\ss}, S.~Peitz, and C.~Sch{\"{u}}tte.
\newblock {Tensor-based dynamic mode decomposition}.
\newblock {\em Nonlinearity}, 31(7):3359--3380, 2018.

\bibitem{WKR15}
M.~O. Williams, I.~G. Kevrekidis, and C.~W. Rowley.
\newblock {A Data-Driven Approximation of the Koopman Operator: Extending
  Dynamic Mode Decomposition}.
\newblock {\em Journal of Nonlinear Science}, 25(6):1307--1346, 2015.

\bibitem{KKS16}
S.~Klus, P.~Koltai, and C.~Sch{\"{u}}tte.
\newblock {On the numerical approximation of the Perron-Frobenius and Koopman
  operator}.
\newblock {\em Journal of Computational Dynamics}, 3(1):51--79, 2016.

\bibitem{PBK15}
J.~L. Proctor, S.~L. Brunton, and J.~N. Kutz.
\newblock {Dynamic mode decomposition with control}.
\newblock {\em SIAM Journal on Applied Dynamical Systems}, 15(1):142--161,
  2015.

\bibitem{PBK18}
J.~L. Proctor, S.~L. Brunton, and J.~N. Kutz.
\newblock {Generalizing Koopman Theory to allow for inputs and control}.
\newblock {\em SIAM Journal on Applied Dynamical Systems}, 17(1):909--930,
  2018.

\bibitem{BBPK16}
S.~L. Brunton, B.~W. Brunton, J.~L. Proctor, and J.~N. Kutz.
\newblock {Koopman invariant subspaces and finite linear representations of
  nonlinear dynamical systems for control}.
\newblock {\em PLoS ONE}, 11(2):1--19, 2016.

\bibitem{KM18b}
M.~Korda and I.~Mezi{\'{c}}.
\newblock {Linear predictors for nonlinear dynamical systems: Koopman operator
  meets model predictive control}.
\newblock {\em Automatica}, 93:149--160, 2018.

\bibitem{KKB17}
E.~Kaiser, J.~N. Kutz, and S.~L. Brunton.
\newblock {Data-driven discovery of Koopman eigenfunctions for control}.
\newblock {\em arXiv:1707.0114}, 2017.

\bibitem{KM18a}
M.~Korda and I.~Mezi{\'{c}}.
\newblock {On Convergence of Extended Dynamic Mode Decomposition to the Koopman
  Operator}.
\newblock {\em Journal of Nonlinear Science}, 28(2):687--710, 2018.

\bibitem{Sag09}
S.~Sager.
\newblock {Reformulations and algorithms for the optimization of switching
  decisions in nonlinear optimal control}.
\newblock {\em Journal of Process Control}, 19(8):1238--1247, 2009.

\bibitem{DPMN18}
N.~Deng, L.~Pastur, M.~Morzy{\'{n}}ski, and B.~R. Noack.
\newblock {Route to chaos in the fluidic pinball}.
\newblock In {\em ASME 2018 Fluids Engineering Division Summer Meeting}, 2018.

\bibitem{Mez13}
I.~Mezi{\'{c}}.
\newblock {Analysis of Fluid Flows via Spectral Properties of the Koopman
  Operator}.
\newblock {\em Annual Review of Fluid Mechanics}, 45:357--378, 2013.

\bibitem{KNK+18}
S.~Klus, F.~N{\"u}ske, P.~Koltai, H.~Wu, I.~Kevrekidis, C.~Sch{\"u}tte, and
  F.~No{\'e}.
\newblock Data-driven model reduction and transfer operator approximation.
\newblock {\em Journal of Nonlinear Science}, 28(3):985--1010, 2018.

\bibitem{AM17}
H.~Arbabi and I.~Mezić.
\newblock {Ergodic Theory, Dynamic Mode Decomposition, and Computation of
  Spectral Properties of the Koopman Operator}.
\newblock {\em SIAM Journal on Applied Dynamical Systems}, 16(4):2096--2126,
  2017.

\bibitem{BD15}
R.~E. Bellmann and E.~D. Stuart.
\newblock {\em {Applied dynamic programming}}.
\newblock Princeton University Press, 2015.

\bibitem{XA00}
X.~Xu and P.~Antsaklis.
\newblock {A dynamic programming approach for optimal control of switched
  systems}.
\newblock In {\em Proceedings of the 39th IEEE Conference on Decision and
  Control}, pages 1822--1827, 2000.

\bibitem{Ell09}
D.~L. Elliott.
\newblock {\em {Bilinear Control Systems - Matrices in Action}}.
\newblock Springer Science + Business Media, 2009.

\bibitem{PP08}
P.~M. Pardalos and V.~Yatsenko.
\newblock {\em {Optimization and Control of Bilinear Systems}}.
\newblock Springer, 2008.

\bibitem{FP02}
J.~H. Ferziger and M.~Peric.
\newblock {\em {Computational Methods for Fluid Dynamics}}.
\newblock Springer Berlin Heidelberg, 3 edition, 2002.

\bibitem{JJT07}
H.~Jasak, A.~Jemcov, and Z.~Tukovic.
\newblock {OpenFOAM : A C++ Library for Complex Physics Simulations}.
\newblock {\em International Workshop on Coupled Methods in Numerical
  Dynamics}, pages 1--20, 2007.

\bibitem{BBP+17}
S.~L. Brunton, B.~W. Brunton, J.~L. Proctor, E.~Kaiser, and J.~N. Kutz.
\newblock {Chaos as an intermittently forced linear system}.
\newblock {\em Nature Communications}, 8(19):1--9, 2017.

\bibitem{AHKO12}
F.~Albrecht, B.~Haasdonk, S.~Kaulmann, and M.~Ohlberger.
\newblock {The localized reduced basis multiscale method}.
\newblock In {\em Proceedings of ALGORITHMY 2012}, pages 393--403, 2012.

\bibitem{LDBK17}
Q.~Li, F.~Dietrich, E.~M. Bollt, and I.~G. Kevrekidis.
\newblock Extended dynamic mode decomposition with dictionary learning: A
  data-driven adaptive spectral decomposition of the koopman operator.
\newblock {\em Chaos: An Interdisciplinary Journal of Nonlinear Science},
  27(10):103111, 2017.

\end{thebibliography}
\end{document}